\documentclass[aos]{imsart} 

\RequirePackage[authoryear]{natbib}
\RequirePackage[colorlinks,citecolor=blue,urlcolor=blue]{hyperref}
\RequirePackage{graphicx}
\usepackage[reqno]{amsmath}
\RequirePackage{amsfonts,amssymb,amsthm, mathtools}
\makeatletter
\usepackage{xcolor}
\usepackage[final]{pdfpages}
\usepackage{ragged2e}
\usepackage{array}
\usepackage{tabularx}
\usepackage{multirow}
\usepackage{rotating}
\usepackage{eurosym}
\usepackage{dsfont}
\usepackage{lmodern}
\usepackage{ulem}
\usepackage{multicol}
\usepackage{hyperref}
\usepackage{etoolbox}
\usepackage{textcomp}
\usepackage{tikz}
\usepackage{tikz-cd}
\usetikzlibrary {positioning}
\usepackage{shadethm}
\usepackage{pdfpages}
\usepackage[utf8]{inputenc}
\usepackage{csquotes}
\usepackage[T1]{fontenc}
\usepackage[german, english]{babel}
\usepackage{latexsym}
\usepackage{nccmath}
\usepackage{color}
\usepackage{comment}
\usepackage{mathrsfs}
\usepackage{clipboard}
\usepackage{autobreak}
\usepackage{hyperref}
\usepackage{bm}
\usepackage{stackengine,scalerel}
\usepackage{cleveref}
\usepackage{stmaryrd}

\makeatother

\renewcommand{\thetable}{\arabic{section}.\arabic{table}}
\numberwithin{equation}{section}

\theoremstyle{plain} 
\newtheorem{theorem}{Theorem}[section]
\newtheorem{lemma}[theorem]{Lemma}
\newtheorem{proposition}[theorem]{Proposition}

\newshadetheorem{assumptionwn}{Assumption A}

\theoremstyle{remark}
\newtheorem{definition}[theorem]{Definition}

\newtheorem{example}[theorem]{Example}

\newtheorem{remark}[theorem]{Remark}

\def\tablename{\footnotesize Table}

\def\ne{\text{ne}}
\def\dis{\text{dis}}
\def\ch{\text{ch}}

\newcommand{\C}{\mathbb{C}} 
\newcommand{\R}{\mathbb{R}} 
\newcommand{\Z}{\mathbb{Z}} 
\newcommand{\N}{\mathbb{N}} 
\newcommand{\BE}{\mathbb{E}}

\newcommand{\BY}{\mathbf{Y}}

\newcommand{\CX}{\mathcal{X}}

\newcommand{\CY}{\mathcal{Y}}
\newcommand{\CZ}{\mathcal{Z}}
\newcommand{\BA}{\mathbf{A}}
\newcommand{\BB}{\mathbf{B}}
\newcommand{\BFC}{\mathbf{C}}

\newcommand{\BS}{\Sigma}

\newcommand{\inde}{\perp\mkern-9.5mu\perp}
\newcommand{\msep}{\bowtie_m}
\newcommand{\bow}{\bowtie}
\newcommand{\beao}{\begin{eqnarray*}}
\newcommand{\eeao}{\end{eqnarray*}\noindent}

\newcommand{\uint}{\int_{-\infty}^{\infty}}

\newcommand{\inst}{\:\raisebox{2pt}{\tikz{\draw[-,densely dashed,line width = 0.5 pt](0,0) -- (5mm,0);}}\:}
\newcommand{\edge}{\:\raisebox{2pt}{\tikz{\draw[-,solid,line width = 0.5pt](0,0) -- (5mm,0);}}\:}

\newcommand{\rarrow}{\:\raisebox{0pt}{\tikz{\draw[->,solid,line width = 0.5 pt](0,0) -- (5mm,0);}}\:}
\newcommand{\larrow}{\:\raisebox{0pt}{\tikz{\draw[->,solid,line width = 0.5 pt](5mm,0) -- (0,0);}}\:}

\newcommand{\ova}{\varepsilon_{a \vert V \setminus \{a,b\}}} 
\newcommand{\ovb}{\varepsilon_{b \vert V \setminus \{a,b\}}}
\newcommand{\ovA}{\varepsilon_{A \vert C}}

\newcommand{\ovB}{\varepsilon_{B \vert C}}
\newcommand{\ovAB}{\varepsilon_{A\cup B \vert C}}
\newcommand{\cab}{c_{\ova \ovb}} 

\newcommand{\cAB}{c_{\ovA \ovB}}
\newcommand{\cAA}{c_{\ovA \ovA}}

\newcommand{\fab}{f_{\ova \ovb}} 
\newcommand{\faa}{f_{\ova \ova}}
\newcommand{\fba}{f_{\ovb \ova}}
\newcommand{\fbb}{f_{\ovb \ovb}}
\newcommand{\fAA}{f_{\ovA \ovA}}
\newcommand{\fAB}{f_{\ovA \ovB}}

\newcommand{\fBB}{f_{\ovB \ovB}}
\newcommand{\fABAB }{f_{\ovAB \ovAB}}
\newcommand{\Rab}{R_{\ova \ovb}} 

\newcommand{\RAB}{R_{\ovA \ovB}}

\allowdisplaybreaks 
\definecolor{darkgreen}{RGB}{0,139,0}
\newcommand{\LS}[1]{{\color{purple} #1}}
\newcommand{\VF}[1]{{\color{darkgreen} #1}}

\begin{document}
\begin{frontmatter}
\title{Partial correlation graphs for \vspace*{0.2cm} \\  continuous-parameter time series}
\runtitle{Partial correlation graphs}

\begin{aug}
{   \author{\fnms{Vicky} \snm{Fasen-Hartmann}\ead[label=e1]{vicky.fasen@kit.edu}\orcid{0000-0002-5758-1999}}
   \and
    \author{\fnms{Lea} \snm{Schenk}\ead[label=e2]{lea.schenk@kit.edu}\orcid{0009-0009-6682-6597}}
}
 \address{Institute of Stochastics, Karlsruhe Institute of Technology\\[2mm] \printead[presep={\ }]{e1,e2}}





\runauthor{V. Fasen-Hartmann and L. Schenk}
\end{aug}

\begin{abstract}
In this paper, we establish the partial correlation graph for multivariate continuous-time stochastic processes, assuming only that the underlying process is stationary and mean-square continuous with expectation zero and spectral density function. In the partial correlation graph,  the vertices are the components of the process and the undirected edges represent partial correlations between the vertices. To define this graph, we therefore first introduce the partial correlation relation for continuous-time processes and provide several equivalent characterisations. In particular, we establish that the partial correlation relation defines a graphoid. The partial correlation graph additionally satisfies the usual Markov properties and the edges can be determined very easily via the inverse of the spectral density function.
 Throughout the paper, we compare and relate the partial correlation graph to the mixed (local) causality graph of \cite{VF23pre}. Finally, as an example, we explicitly characterise and interpret the edges in the partial correlation graph for the popular multivariate continuous-time AR (MCAR) processes.
\end{abstract}



\end{frontmatter}

\section{Introduction}\label{sec: Introduction}

Our interest in this paper is in graphical models for wide-sense stationary and mean-square continuous stochastic processes. Graphical models are probabilistic networks, where the vertices represent the components of a random object, e.g., a random vector or a vector-valued stochastic process, and the edges illustrate specific interconnections between them. They are popular because they visualise dependency structures of the random object in a clear and simple way,  which can then be analysed, interpreted, and easily communicated. Furthermore, graphical models are an important tool for dimension reduction in high-dimensional models. Due to the growth of complex multivariate data sets and networks, the theory and methodology of graphical models have experienced a surge of research development in probability theory and statistics \citep{WI08,Edwards,LA04,Handbook:graphical}, and they have been applied in fields as diverse as biology, neuroscience, economics, finance, and psychology, to name just a few.

Although in networks of interconnected processes the data are observed in discrete time, in many situations it is more appropriate to specify the underlying stochastic process in \textit{continuous time}. This is particularly necessary for high-frequency data, irregularly spaced data or data with missing observations, which are common in finance, econometrics, signal processing, and turbulence. In addition, many physical and signal processing models are formulated in continuous time, so such an approach is often more natural.

Overall, however, there is very little theory on graphical models for multivariate stochastic processes in continuous time. The established graphical models are mostly limited to conditional independence and local independence graphs, which have been studied by \cite{Mogensen:Hansen:2020, Mogensen:Hansen:2022, DI07, DI08, Aalen87, Schweder}. They are particularly suitable for semimartingales and point processes, but do not seem to be the right tool for time series. A general approach for graphical continuous-parameter time series models are the (local) causality graphs of \cite{VF23pre, VF23preb}, which are mixed graphs representing Granger causalities and contemporaneous correlations. The causality graphs satisfy the usual Markov properties and the theory holds for a very large class of time series models.  However, the computation of the edges can be quite challenging in certain examples and the characterisations may not be convenient, as we know for multivariate continuous-time ARMA processes from \cite{VF23preb}, which is problematic in practice. Until now, an undirected graphical model for continuous-time stationary processes has been lacking. Therefore, we aim to fill this gap and provide a user-friendly and powerful undirected graphical model for continuous-time stochastic processes.

In the graphical model we propose in this paper, the vertices $V=\{1,\ldots,k\}$ represent the components of a $k$-dimensional process and the edges visualise partial correlations between these components. The concept of \textit{partial correlation} is an important and well-studied measure of dependence in statistics. For an $\R^k$-valued random vector \mbox{$\BY=(Y_1,\ldots,Y_k)^\top$} with positive definite covariance matrix $\Sigma$, the partial correlation of $Y_a$ and $Y_b$ given  $\overline\BY\coloneqq(Y_l)_{l\in\{1,\ldots,k\}\backslash\{a,b\}}$ measures the correlation of the real-valued random variables $Y_a$ and $Y_b$ after removing the linear effects of $\overline\BY$. The partial correlation is determined as follows: Denote by $\Sigma_{AB}$ the respective submatrix of $\Sigma$ for $A, B\subseteq V$ and consider the linear regression problems
\begin{align} \label{regression}
    \beta_l&=\text{argmin}_{\beta\in\R^{k-2}}\BE (Y_l-\beta^\top\overline \BY)^2, \quad l\in\{a,b\}.
\intertext{These problems have the well-known solution \citep{FUS:book:2010,Anderson} }
    \beta_l&= \left(\Sigma_{V\setminus\{a,b\} V\backslash\{a,b\}} \right)^{-1}\Sigma_{V\backslash\{a,b\} l}.  \nonumber
\intertext{Furthermore, the residuals
    $\varepsilon_{a|V\backslash\{a,b\}} \coloneqq Y_a-\beta_a^\top\overline \BY$  and $\varepsilon_{b|V\backslash\{a,b\}} \coloneqq Y_b-\beta_b^\top\overline \BY$
 satisfy }
    \text{Cov}(\varepsilon_{a|V\backslash\{a,b\}},\varepsilon_{b|V\backslash\{a,b\}})&=\Sigma_{ab}-\Sigma_{aV\backslash\{a,b\}} \left(\Sigma_{V\backslash\{a,b\}V\backslash\{a,b\}} \right)^{-1}
    \Sigma_{V\backslash\{a,b\} b}, \label{1.2}
\end{align}
which is the\textit{ partial covariance }of $Y_a$ and $Y_b$ given $Y_{V\setminus \{a,b\}}$. Similarly, the correlation of the residuals is called \textit{partial correlation} of $Y_a$ and $Y_b$ given $ Y_{V\setminus \{a,b\}}$, also known as \textit{coherence}, and is equal to
\begin{align} \label{1.3}
    \text{Corr}( \varepsilon_{a|V\backslash\{a,b\}}, \varepsilon_{b|V\backslash\{a,b\}})
    &= \frac{\text{Cov}( \varepsilon_{a|V\backslash\{a,b\}}, \varepsilon_{b|V\backslash\{a,b\}})}{\sqrt{\text{Cov}( \varepsilon_{a|V\backslash\{a,b\}}, \varepsilon_{a|V\backslash\{a,b\}})\text{Cov}( \varepsilon_{b|V\backslash\{a,b\}}, \varepsilon_{b|V\backslash\{a,b\}})}} \nonumber \\
    &= -\frac{[\Sigma^{-1}]_{ab}}{\sqrt{[\Sigma^{-1}]_{aa}[\Sigma^{-1}]_{bb}}}.
\end{align}
From the representation \eqref{1.3} we see that the partial correlation is completely determined by the \textit{concentration }(precision) \textit{matrix} $\Sigma^{-1}$. For a Gaussian random vector, zero partial correlation is even equivalent to  $Y_a$ and $Y_b$ being independent given  $\overline\BY$.

An extension of partial correlation to stationary time series models in \textit{discrete time} is quite old \citep{Tick} and is ubiquitous in the analysis of multivariate time series  \citep{Priestley, BR01, Gardner:1988}. Recall that wide-sense stationary processes are those where the process has a constant expectation at each time point and the existing covariance function depends only on the time lags. Moreover, the spectral representation of a multivariate time series is a frequency domain representation and the spectral density is directly related to the autocovariance function in the time domain by Fourier transformation.
For these time series models, the partial covariance function is zero if and only if the partial spectral density function is zero, such that in the frequency domain, the role of the partial correlation function is taken over by the \textit{spectral coherence function} of the noise process, the normalised cross spectral density, and the role of the covariance matrix $\Sigma$  is taken over by the matrix-valued spectral density function of the process.
Here, in the context of time series, the spectral coherence function of the noise process measures the linear dependence between two components of a multivariate time series after removing the linear effects of the remaining components in the frequency domain, and in \Cref{remainder satisfies Assumption 1} we present the corresponding result to \eqref{1.2}  and \Cref{characterisation via inverse density} corresponds to \eqref{1.3}, respectively. The applications of spectral coherence are very broad, especially in signal processing, but the word coherence may have a slightly different meaning in different fields  \citep{GARDNER1992113}. However, to the best of our knowledge, a mathematically rigorous theory for the definition of partial correlation for continuous-parameter time series is missing in the literature, so we include the theory first and relate it to an optimisation problem as in \eqref{regression} in \Cref{Loesung ist Projektion} and \Cref{Optimisation problem Dahlhaus}. 
It is important to note that in the above regression problem, $\beta_l^\top\overline \BY$ is the linear projection of $Y_l$ on the linear space generated by the components of $\overline{\BY}$ because our approach builds on this idea. In particular, we show that our definition of partial correlation satisfies the important graphoid properties.

The subject of this paper are then \textit{partial correlation graphs} for continuous-time wide-sense stationary and mean-square continuous stochastic processes with expectation zero and spectral density.
Partial correlation graphs for discrete-time wide-sense stationary stochastic processes with expectation zero and spectral density originated in \cite{BR96}  and \cite{DA00} and are a widely used frequency domain approach for constructing graphs. In our  graphical model and in the model of \cite{DA00}, the vertices are the components of a multivariate time series  and the edges between the vertices are drawn when the
spectral coherence function in these components is the zero function, meaning that the component processes are partially uncorrelated given the remaining process.
 The method of \cite{DA00} has since been used in a wide variety of applications, including the identification of synaptic connections in  air pollution data \citep{DA00}, human tremor data \citep{DA03}, vital signs of intensive care patients \citep{GA00}, financial data \citep{Abdelwahed}, and neuro-physical signals 
 \citep{DA97,EI03,MEDKOUR2009374}, which demonstrates the popularity of partial correlation graphs in identifying a network structure.

This paper aims to define a probabilistic network of interconnected continuous-time stochastic processes, where the dependence structure in the network is modelled by partial correlation. The proposed partial correlation graph is simple in the sense that there are neither loops from a vertex to itself nor any multiple edges between vertices and it satisfies the required Markov properties that associate the graph factorisation to the partial correlation. Moreover, it is easy to handle in applications because the edges reflect zero entries in the inverse spectral density function. We derive important relations between the undirected partial correlation graph and the recently introduced mixed causality graph of \cite{VF23pre}.
In the mixed causality graph, the directed and undirected edges can be defined by conditional orthogonality relations of properly defined linear subspaces generated by the underlying stochastic process, similarly, the edges in the partial correlation graph can be defined by conditional orthogonality. We use this commonality to compare both graphical models and to show the important connection that the edges in the partial correlation graph are also edges in the augmented causality graph. Furthermore, as an example, we apply the partial correlation graph to multivariate continuous-time autoregressive (MCAR) processes and present a perspective on estimation. In the context of MCAR processes, we additionally obtain that the edges of the partial correlation graph are also edges in the corresponding augmented local causality graph of \cite{VF23pre}. Finally, a major conclusion of this paper is that the edges of the continuous-time model are in general not identifiable from equidistantly sampled observations, but this is different for high-frequency data.

\subsection*{Structure of the paper}
The paper is structured as follows. In \Cref{sec: Preliminaries}, we lay the groundwork for the paper by introducing relevant properties of multivariate wide-sense stationary and mean-square continuous processes. Then, in \Cref{sec: Partial correlation}, we define the partial correlation relation and establish characterisations and properties. This preliminary work results in the definition of the partial correlation graph \mbox{$G_{PC}=(V, E_{PC})$} in \Cref{sec: Partial correlation graphs}, where we also discuss edge characterisations, as well Markov properties, and the relations to the causality graph. As an example, in \Cref{sec: Partial correlation graphs for MCAR processes}, we apply the partial correlation graph to MCAR processes and compare it to the local causality graph for MCAR processes. Finally, we complete the paper with a brief conclusion in \Cref{Conclusion}. The proofs of the paper are given in \Cref{sec:proofs}.

\subsection*{Notation}
In the following, $I_{k}\in \R^{k\times k}$ is the $(k\times k)$-dimensional identity matrix, $0_{k\times d} \in \R^{k\times d}$ is the $(k\times d)$-dimensional zero matrix, and $0_{k}$ is either the $k$-dimensional zero vector or the $(k\times k)$-dimensional zero matrix, which should be clear from the context. The vector $e_a\in \R^k$ is the $a$-th unit vector. The entries and submatrices of a matrix $M$ are denoted by $\left[M\right]_{ab}$ for $a,b\in V$ and $\left[M\right]_{AB}$ for $A,B\subseteq V$, respectively. The cardinality of a set $A$ is denoted by $|A|$. For hermitian matrices $M, N \in \C^{k\times k}$, we write $M \geq N$ if and only if $M-N$ is positive semi-definite. Similarly, we write $M>0$ if and only if $M$ is positive definite and define $\sigma(M)$ as the set of eigenvalues of $M$. Finally, for a function $f: \R \rightarrow \C^{k\times k}$ with $f(\lambda)>0$, we define the function $g: \R \rightarrow \C^{k\times k}$ by $g(\lambda)= f(\lambda)^{-1}$, $\lambda \in \R$.

\section{Preliminaries}\label{sec: Preliminaries}
In this paper, we consider wide-sense stationary and mean-square continuous stochastic processes $\CY_V=\left(Y_V\left(t\right)\right)_{t \in \R}$ in continuous time with index set $V=\{1,\ldots,k\}$, $\BE(Y_V(t))=0_k$, and a spectral density function $f_{Y_VY_V}(\lambda)$ for $\lambda \in \R$. Note that $\CY_V$ is mean-square continuous if and only if
    \begin{align}\label{symmetrie of covariance}
    \lim_{t \rightarrow 0} c_{Y_VY_V}(t) = c_{Y_VY_V}(0).
    \end{align}
The autocovariance function of $\CY_V$ is denoted by $(c_{Y_VY_V}(t))_{t\in\R}=(\BE[Y_V(t)Y_V(0)^\top])_{t\in\R}$.

In this section, we present well-known properties of these processes that are relevant to this work. The results date back to \cite{KH34} and \cite{CR39} and were summarised in a comprehensive overview, e.g., by \cite{DO60} and \cite{RO67}.

A key property of wide-sense stationary and mean-square continuous stochastic processes with expectation zero and existing spectral density function is their spectral representation
    \begin{align}\label{spectral representation of stationary process}
    Y_V(t) = \uint e^{i \lambda t} \Phi_V(d\lambda), \quad t\in \R,
    \end{align}
with respect to a random orthogonal measure $\Phi_V=(\Phi_1,\ldots,\Phi_k)^\top$, where
    \begin{align*}
    \BE[\Phi_V(d\lambda) \overline{\Phi_V(d\mu)}^\top] = \delta_{\lambda=\mu} f_{Y_VY_V}(\lambda) d\lambda, \quad
    \BE[\Phi_V(d\lambda)] = 0_k,
    \end{align*}
and $\delta_{\lambda=\mu}$ is the Kronecker Delta. We refer to the function $f_{Y_AY_B}(\lambda)=[f_{Y_VY_Y}(\lambda)]_{AB}$, \linebreak $\lambda\in\R$, as the {cross-spectral density function} of the subprocesses $\CY_A$ and $\CY_B$ with \linebreak $A,B\subseteq V$. Important properties of the spectral density function are the following.


\begin{lemma}\label{properties of spectral density}
Let   $\lambda, t \in \R$. Then the following statements hold.
\begin{itemize}
    \item[(a)] $\uint \| f_{Y_VY_V}(\lambda) \|\lambda <\infty$,
    \item[(b)] $c_{Y_VY_V}(t) = \uint e^{i\lambda t} f_{Y_VY_V}(\lambda) d\lambda$,
    \item[(c)] $f_{Y_VY_V}(\lambda)\geq 0$ and $f_{Y_VY_V}(\lambda) = \overline{f_{Y_VY_V}(\lambda)}^\top$.
\end{itemize}
\end{lemma}

\begin{remark}
To obtain the one-to-one relationship
\begin{align}\label{Sufficient condition for Assumption 3}
    f_{Y_VY_V}(\lambda) = \frac{1}{2\pi} \int_{-\infty}^{\infty}  e^{-i \lambda t} c_{Y_VY_V}(t) dt, \quad \lambda \in \R,
\end{align}
between $c_{Y_VY_V}(t)$ and $f_{Y_VY_V}(\lambda)$ by Fourier transformation, additional integrability assumptions on the covariance function are required.
For example, suppose that 
$\int_{-\infty}^{\infty} \| c_{Y_VY_V}(t) \| dt < \infty$, then $\CY_V$ has a spectral density function given by \eqref{Sufficient condition for Assumption 3}. However, in this paper, we only require the existence of a spectral density function and not the relation \eqref{Sufficient condition for Assumption 3}. Therefore,  long memory processes such as multivariate fractionally integrated CARMA processes are also covered in this paper.
\end{remark}

In addition to the spectral density function, we introduce the spectral coherence function of $\CY_A$ and $\CY_B$, which is obtained by rescaling the cross-spectral density function of $\CY_A$ and $\CY_B$, and provides a measure of the strength of the dependence. 

\begin{definition}\label{partial   correlation matrix}
The \textsl{spectral coherence function} of $\CY_A$ and $\CY_B$ is defined as
\begin{align*}
R_{Y_A Y_B}(\lambda) \coloneqq \Bigl( f_{Y_A Y_A}(\lambda) \Bigr)^{-1/2} f_{Y_A Y_B}(\lambda) \Bigl(  f_{Y_B Y_B}(\lambda) \Bigr)^{-1/2}, \quad \lambda \in \R.
\end{align*}
If $f_{Y_A Y_A}(\lambda)$ or $f_{Y_B Y_B}(\lambda)$ is singular for some $\lambda \in \R$, we set $R_{Y_A Y_B}(\lambda) \coloneqq 0_{|A| \times |B|}$.
\end{definition}

The following linear spaces generated by subprocesses
 $\CY_C$ of $\CY_V$,  $C\subseteq V$, are also regularly used in this paper. For $t\in \R$, we define
\begin{align*}
L_{Y_C}(t) &\coloneqq \overline{ \left\{ \sum_{c\in C}\gamma_c Y_c(t): \gamma_c \in \C \right\} },\\
\mathcal{L}_{Y_C}(t) &\coloneqq  \overline{ \left\{ \sum_{i=1}^\ell \sum_{c\in C} \gamma_{c,i} Y_c(t_i): \gamma_{c,i} \in \C, \: -\infty < t_1 \leq \ldots \leq t_\ell \leq t, \: \ell \in \N \right\} }, \\
\mathcal{L}_{Y_C} &\coloneqq \overline{ \left\{ \sum_{i=1}^\ell \sum_{c\in C} \gamma_{c,i} Y_c(t_i): \gamma_{c,i} \in \C, \: -\infty < t_1 \leq \ldots \leq t_\ell <\infty, \: \ell \in \N \right\} },
\end{align*}
where $\overline{ \{\cdot\} }$ denotes the mean-square closure. Moreover, let $\Phi_C$ be the random spectral measure from the spectral representation \eqref{spectral representation of stationary process} of $\CY_C$ and
\begin{align*}
    L^2\left(f_{Y_C Y_C}\right) \coloneqq \Bigg\{\varphi^\top : \R \rightarrow \C^{|C|}: \varphi \text{ measurable,} \uint \left\vert \varphi(\lambda) f_{Y_C Y_C}(\lambda) \overline{\varphi(\lambda)}^\top \right\vert d\lambda < \infty \Bigg\}.
\end{align*}
Then we define as well the space
\begin{align*}
    \mathcal{L}_{Y_C}^* &\coloneqq \left\{ \uint e^{i \lambda t} \varphi(\lambda) \Phi_C(d\lambda) : \varphi \in L^2\left(f_{Y_C Y_C}\right), \: t\in \R \right\},
\end{align*}
which is related to $\mathcal{L}_{Y_C}$ as follows.


\begin{lemma}\label{Gleichheit linearer Raeume}
The spaces $\mathcal{L}_{Y_C}$ and $\mathcal{L}_{Y_C}^*$ are equal.
\end{lemma}
For the considerations of the partial correlation relation, this equivalence is crucial.
\section{Partial correlation relation}\label{sec: Partial correlation}
Next, we introduce the concept of partial correlation for wide-sense stationary, mean-square continuous processes $\CY_V$ that have expectation zero and a spectral density function. 
Therefore, in \Cref{Partial correlation and conditional orthogonality}, we define and interpret the partial correlation relation and compute the orthogonal projections therein. Additionally, we discuss properties of $Y_A(t)$ given the linear information of $\CY_C$, i.e., the resulting noise process. \Cref{Section: Characterisations of partial correlation} is then devoted to characterisations of the partial correlation relation. We provide characterisations in terms of the spectral density function and the spectral coherence function of the noise processes. Importantly, we present the key characterisation involving the inverse of the spectral density function of the underlying process $\CY_V$. We conclude the section with the main result of this section that the partial correlation relation satisfies the graphoid properties. Throughout this section, $A, B, C$ are subsets of $V$.

\subsection{Partial correlation relation and orthogonal projections}\label{Partial correlation and conditional orthogonality}
Let us introduce the concept of partial correlation and make some comments on that definition.

\begin{definition}\label{Definition of partial correlation}
The two subprocesses  $\CY_A$ and $\CY_B$ of $\CY_V$ are defined to be \textsl{partially uncorrelated} given $\CY_C$ if and only if
\begin{align*}
    \BE \left[ \left(  Y_a(t)-P_{\mathcal{L}_{Y_C}}  Y_a(t) \right) \overline{\left(  Y_b(t)-P_{\mathcal{L}_{Y_C}}  Y_b(t)\right)} \right] = 0 \quad \forall \:  a\in A, \: b\in B, \: t\in \R,

\end{align*}
where $P_{\mathcal{L}_{Y_C}}$ is the orthogonal projection on $\mathcal{L}_{Y_C}$. In short, we write $\CY_A \inde  \CY_B \: \vert \:  \CY_{C}$.
\end{definition}

\begin{remark}\mbox{}\label{Bedingte Orthogonalitaet}
    \item[(a)] The partial uncorrelation $\CY_A \inde \CY_B \: \vert \:  \CY_{C}$ states, as desired, that, for all $t\in \R$, $Y_A(t)$ and $Y_B(t)$ are uncorrelated given the linear information provided by $\CY_C$ over all time points. The concept can be seen as an extension of the definition of partial correlation for random vectors in \Cref{sec: Introduction}. In terms of the conditional orthogonality relation $\perp$ (cf.~\citealp{EI07}, Appendix A), this means that
\begin{align*} 
\CY_A \inde  \CY_B \: \vert \:  \CY_{C}
\quad &\Leftrightarrow \quad
\BE \left[ \left( Y^A-P_{\mathcal{L}_{Y_C}} Y^A \right) \left(  Y^B-P_{\mathcal{L}_{Y_C}} Y^B \right)\right] = 0 \nonumber \\
    \quad & \hspace{1,1cm} \forall \:  Y^A \in L_{Y_A}(t), \; Y^B \in L_{Y_B}(t), \: t\in \R, \nonumber \\
    \quad & \Leftrightarrow \quad
    L_{Y_A}(t) \perp L_{Y_B}(t) \: \vert \: \mathcal{L}_{Y_C} \quad \forall \: t\in \R.
\end{align*}
    \item[(b)] Certainly, the partial correlation relation is symmetric and
\begin{align*}
    \CY_A \inde  \CY_B \: \vert \:  \CY_{C}
    \quad \Leftrightarrow \quad
    \CY_a \inde  \CY_b \: \vert \:  \CY_{C},   \quad \forall \: a\in A, b\in B,
\end{align*}
which is useful for verifying zero partial correlation. Furthermore, statements can usually be made without loss of generality for $A=\{a\}$ and $B=\{b\}$. The corresponding results in the multivariate case follow immediately.
\end{remark}

In order to work with the partial correlation relation, we compute the orthogonal projections $P_{\mathcal{L}_C} Y_a(t)$ in the next proposition. Therefore, remark that stochastic integrals of deterministic Lebesgue measurable functions with respect to a random orthogonal measure are defined in the usual $L^2$-sense. For details on the definition and properties of such integrals, we refer to  \cite{DO60} and \cite{RO67}.


\begin{proposition}\label{Loesung ist Projektion}
Suppose that $f_{Y_CY_C}(\lambda)>0$ for $\lambda \in \R$. Then, for $t\in \R$, the orthogonal projection is equal to
    \begin{align*}
        P_{\mathcal{L}_C} Y_a(t)
        = \uint e^{i \lambda t} f_{Y_a Y_C}(\lambda) f_{Y_CY_C}(\lambda)^{-1} \Phi_C(d \lambda),
    \end{align*}
where $\Phi_C$ is the random spectral measure from the spectral representation \eqref{spectral representation of stationary process} of $\CY_C$. Furthermore,  $P_{\mathcal{L}_C} Y_a(t)$ is the solution to the optimisation problem
   \begin{align}\label{optimisation problem}
    \underset{\varphi_{a \vert C} \in L^2 \left(f_{Y_C Y_C} \right)}{\text{min}}
    \BE \left[ \left\vert Y_a(t) - \uint e^{i \lambda t} \varphi_{a \vert C}(\lambda) \Phi_C(d \lambda) \right\vert^2
    \right].
    \end{align}
    Finally,
 $P_{\mathcal{L}_C} Y_A(t)=(P_{\mathcal{L}_C} Y_a(t))_{a\in A}$ can be calculated component-wise.

\end{proposition}
Note that the requirement for the existence of a partially positive definite spectral density function allows for an explicit representation of the orthogonal projection.

\begin{remark}\label{Optimisation problem Dahlhaus}
The choice of the term partial correlation relation is inspired by the partial correlation relation for discrete-time stationary processes in \cite{BR01} and  \cite{DA00}. However, the discrete-time concept is motivated by an optimisation problem similarly to \eqref{regression} (\citealp{BR01}, Theorem 8.3.1 and  \citealp{DA00}, relation (2.1) and Definition 2.1). To see the correspondence, suppose that the function $\varphi_{a \vert C}$ in the optimisation problem \eqref{optimisation problem} is the Fourier transform of an \textsl{integrable} function $d_{a \vert C}(t)$, $t\in \R$. Then, for $t\in \R$, \cite{RO67}, I, Example 8.3, provides
\begin{align*}
    \uint e^{i \lambda t} \varphi_{a \vert C}(\lambda) \Phi_C(d \lambda)
    = \uint d_{a \vert C}(t-s) Y_C(s) ds.
\end{align*}
With this representation, we have the similarity of our optimisation problem \eqref{optimisation problem} to the discrete-time optimisation problem
\begin{align*}
    \underset{d_{a \vert C}}{\text{min }} \BE \left[ \left\vert Z_a(t) - \sum_{u=-\infty}^\infty d_{a \vert C}(t-u) Z_{C} (u) \right\vert^2 \right],
\end{align*}
and to \eqref{regression}.
The advantage of our approach is that we require weaker assumptions.
Given this parallelism, similarities with \cite{DA00} are to be expected in various sections of this paper.

\end{remark}

Finally, we define the multivariate noise process
\begin{align*}
\varepsilon_{A \vert C}(t)
        \coloneqq Y_A(t) - P_{\mathcal{L}_C} Y_A(t)
        = Y_A(t) - \uint e^{i \lambda t} f_{Y_A Y_C}(\lambda) f_{Y_CY_C}(\lambda)^{-1} \Phi_C(d \lambda), \quad t\in\R,
\end{align*}
with the following crucial properties.


\begin{lemma}\label{remainder satisfies Assumption 1}
Suppose that $f_{Y_CY_C}(\lambda)>0$ for  $\lambda \in \R$. Then the noise processes $(\ovA(t))_{t\in \R}$ and $(\ovB(t))_{t\in \R}$ are wide-sense stationary and stationary correlated with (cross-) spectral density function
\begin{align*}
\fAB(\lambda)&=f_{Y_A Y_B}(\lambda) - f_{Y_A Y_C}(\lambda) f_{Y_C Y_C}(\lambda)^{-1} f_{Y_C Y_B}(\lambda) \quad \text{for almost all $\lambda \in \R$},
\end{align*}
 and (cross-) covariance function
 \begin{align*}
\cAB(t) = \uint e^{i \lambda t} \left(f_{Y_A Y_B}(\lambda) - f_{Y_A Y_C}(\lambda) f_{Y_CY_C}(\lambda)^{-1}  f_{Y_C Y_B}(\lambda)  \right) d\lambda  \quad \text{for all $t \in \R$}.
\end{align*}
\end{lemma}

\subsection{Characterisations of the partial correlation relation}\label{Section: Characterisations of partial correlation}

In this section, we present several characterisations of the partial correlation relation. We start with simple characterisations in terms of the (cross-) covariance function, the (cross-) spectral density function, and the spectral coherence function of the noise processes, analogous to the discrete-time results in Remark 2.3 of \cite{DA00}.

\begin{proposition}\label{characterisation with spectral density function}
Suppose that $f_{Y_CY_C}(\lambda)>0$ for $\lambda \in \R$. Then the following equivalences hold.
    \begin{align*}
        \CY_A \inde  \CY_B \: \vert \:  \CY_{C}
        \quad &\Leftrightarrow \quad
        \cAB(t) = 0_{|A| \times |B|} \quad \: \text{for all $t \in \R$,} \\
        \quad &\Leftrightarrow \quad
        \fAB(\lambda) = 0_{|A| \times |B|} \quad \text{for almost all $\lambda \in \R$.}
    \end{align*}
In particular, these conditions imply that the spectral coherence function satisfies $\RAB(\lambda) = 0_{|A| \times |B|}$ for almost all $\lambda \in \R$. If $\fAA(\lambda)>0$ and $\fBB(\lambda)>0$ for $\lambda \in \R$, then the converse holds as well.
\end{proposition}

\begin{remark}\mbox{}\label{sufficient criterion for non-singular}
\begin{itemize}
\item[(a)] The assumption that $\fAA(\lambda)>0$ for $\lambda\in \R$ excludes the case \mbox{$\ovA(t)=0_{|A|}$} \linebreak $\mathbb{P}$-a.s.~for $t\in \R$, e.g., the case where $Y_a(t) \in \mathcal{L}_{Y_C}$ for $a\in A$. This can be explained as follows. If  $\ovA(t)=0_{|A|}$ $\mathbb{P}$-a.s.~for $t\in \R$, then 
$\fAA(\lambda)=0_{|A|}\in \R^{|A| \times |A|}$ is not positive definite for $\lambda\in \R$. We can therefore assume that $A\cap C=\emptyset$.
\item[(b)] For $A\cap C=\emptyset$, \cite{BE09}, Proposition 8.2.4 provides  that $f_{Y_{A\cup C} Y_{A\cup C}}(\lambda)>0$ if and only if $f_{Y_{C} Y_{C}}(\lambda)>0$ and $\fAA(\lambda)>0$.
\item[(c)] If $A\cap C=\emptyset$ and $f_{Y_{A \cup C}Y_{A \cup C}}(\lambda)>0$ for $\lambda \in \R$, then $\fAA(\lambda)>0$ and \Cref{characterisation with spectral density function} results in $\CY_A \not \inde  \CY_A \: \vert \:  \CY_{C}$. In the following, we always assume a sufficient condition for $f_{Y_{A\cup C}Y_{A\cup C}}(\lambda)>0$, so we can also exclude the case $A\cap B=\emptyset$ from our analysis and assume throughout the remaining section that $A, B, C \subseteq V$ are disjoint. \qedhere
\end{itemize}
\end{remark}

Finally, we present a very simple characterisation of the partial correlation relation in terms of the inverse of the spectral density function, which we denote, for $A\subseteq V$ and $\lambda \in \R$, by
\begin{align*}
  g_{Y_{A}Y_{A}}(\lambda)
\coloneqq f_{Y_{A}Y_{A}}(\lambda)^{-1}.
\end{align*}
The corresponding discrete-time result is given in Theorem 2.4 of \cite{DA00} and we refer to the proof there.

\begin{proposition}\label{characterisation with inverse}
Suppose that $A,B,C \subseteq V$ are disjoint and \mbox{$f_{Y_{A\cup B \cup C}Y_{A\cup B \cup C}}(\lambda)>0$} for $\lambda\in \R$. Then the following equivalence holds.
\begin{align*}
\CY_A \inde  \CY_B \: \vert \:  \CY_{C}
  & \quad \Leftrightarrow \quad
\left[g_{Y_{A\cup B \cup C}Y_{A\cup B \cup C}}(\lambda)\right]_{AB} = 0_{|A| \times |B|} \quad \text{for almost all $\lambda \in \R$.}
\end{align*}
\end{proposition}

The characterisation via the inverse of the spectral density function of $\CY_{A\cup B \cup C}$ can be used to explain the effect of an unobserved multivariate process $\CY_C$, a so-called confounder process. The following lemma introduces a relationship between the inverse of the spectral density function of a full process $\CY_V$ and the inverse of the spectral density function of a process reduced $\CY_V$ by a confounder process $\CY_C$. This result is the continuous-time counterpart to
\cite{DA00}, Remark 2.5. Since it is a straightforward calculation, we omit the proof. 

\begin{lemma}\label{deleting confounder}
Suppose that $A,B,C \subseteq V$ are disjoint and $f_{Y_{V}Y_{V}}(\lambda)>0$ for $\lambda\in \R$. Then, for $\lambda \in \R$, the equality
\begin{align*}
\big[g_{Y_{V\setminus C}Y_{V\setminus C}}(\lambda)\big]_{AB}
=\left[g_{Y_VY_V}(\lambda)\right]_{AB} - \left[g_{Y_VY_V}(\lambda)\right]_{AC} \left(\left[g_{Y_VY_V}(\lambda)\right]_{CC}\right)^{-1} \left[g_{Y_VY_V}(\lambda)\right]_{CB}
\end{align*}
holds.
\end{lemma}

\begin{remark}
For an interpretation of this result (cf.~Remark 2.5 in \citealp{DA00}), we analyse the case
\begin{align*}
\big[g_{Y_{V\setminus c}Y_{V\setminus c}}(\lambda)\big]_{ab}
=\left[g_{Y_VY_V}(\lambda)\right]_{ab} - \left[g_{Y_VY_V}(\lambda)\right]_{ac} \left(\left[g_{Y_VY_V}(\lambda)\right]_{cc}\right)^{-1} \left[g_{Y_VY_V}(\lambda)\right]_{cb}.
\end{align*}
This equation explains the relation between the partial correlation structure in the full process $\CY_{V}$ and the partial correlation structure in the reduced process $\CY_{V\setminus \{c\}}$: If $\CY_a$ and $\CY_b$ are partially uncorrelated given $\CY_{V\setminus \{a,b\}}$ ($\left[g_{Y_VY_V}(\lambda)\right]_{ab} = 0 $ for almost all $\lambda \in \R$), but there is a partial correlation between $\CY_a$ and $\CY_c$ given $\CY_{V\setminus \{a,c\}}$ and between $\CY_c$ and $\CY_b$ given $\CY_{V\setminus \{b,c\}}$ with $\left[g_{Y_VY_V}(\lambda)\right]_{ac} \neq 0$ and $\left[g_{Y_VY_V}(\lambda)\right]_{cb}\neq 0$ on some non-zero set, this causes a partial correlation between $\CY_a$ and $\CY_b$ given $\CY_{V\setminus \{a,b,c\}}$ ($\big[g_{Y_{V\setminus c}Y_{V\setminus c}}(\lambda)\big]_{ab}\neq 0$) in the reduced process $\CY_{V\setminus \{c\}}$.
\end{remark}

Finally, we establish the main result of this section, namely that the partial correlation relation satisfies the graphoid properties.

\begin{proposition}\label{graphoid}
Suppose that $A,B,C,D \subseteq V$ are disjoint and $f_{Y_VY_V}(\lambda)>0$ for \mbox{$\lambda\in \R$.} Then the partial correlation relation defines a graphiod, i.e., it satisfies the following properties:
\begin{itemize}
  \item[(P1)] \makebox[2.6cm][l]{Symmetry:} $\CY_A \inde \CY_B \: \vert \: \CY_C$ $\Rightarrow$ $\CY_B \inde \CY_A \: \vert \: \CY_C$.
  \item[(P2)] \makebox[2.6cm][l]{Decomposition:} $\CY_A \inde \CY_{B \cup C} \: \vert \: \CY_D$ $\Rightarrow$ $\CY_A \inde \CY_B \: \vert \: \CY_D$.
  \item[(P3)] \makebox[2.6cm][l]{Weak union:} $\CY_A \inde \CY_{B \cup C} \: \vert \: \CY_D$ $\Rightarrow$ $\CY_A \inde \CY_B \: \vert \: \CY_{C \cup D}$.
  \item[(P4)] \makebox[2.6cm][l]{Contraction:} $\CY_A \inde \CY_B \: \vert \: \CY_D$ and $\CY_A \inde \CY_C \: \vert \: \CY_{B \cup D}$ $\Rightarrow$ $\CY_A \inde \CY_{B \cup C} \: \vert \: \CY_D$.
  \item[(P5)] \makebox[2.6cm][l]{Intersection:} $\CY_A \inde \CY_B \: \vert \: \CY_{C \cup D}$ and $\CY_A \inde \CY_C \: \vert \: \CY_{B \cup D}$ $\Rightarrow$ $\CY_A \inde \CY_{B \cup C} \: \vert \: \CY_D$.
\end{itemize}
\end{proposition}

\begin{remark}\mbox{}
\begin{itemize}
    \item[(a)] The property (P1) is immediately clear. The properties (P2), (P3) and (P5) were already established by \cite{DA00} in Lemma 3.1. For (P4) we apply \Cref{characterisation with inverse} and \Cref{deleting confounder} which is done in the appendix.
    \item[(b)] Although the partial correlation can be characterised by conditional orthogonality (\Cref{Bedingte Orthogonalitaet}),  the results of \cite{VF23pre} are not directly applicable. The reason is the following: Due to \Cref{Bedingte Orthogonalitaet}, $\CY_A \inde \CY_{B \cup C} \: \vert \: \CY_D$ is equivalent to the conditional orthogonality relation $L_{Y_A}(t) \perp L_{Y_{B \cup C}} (t) \: \vert \: \mathcal{L}_{Y_D}$ for all $t\in \R$. Thus the weak union property of the conditional orthogonality relation \cite[Lemma 2.2]{VF23pre} gives $L_{Y_A}(t) \perp L_{Y_{B}} (t) \: \vert \: \overline{\mathcal{L}_{Y_D}+ L_{Y_{C}} (t)}$ for all $t\in \R$. This is not the same as $L_{Y_A}(t) \perp L_{Y_{B}} (t) \: \vert \: \mathcal{L}_{Y_{C \cup D}}$ for all $t\in \R$, i.e., $\CY_A \inde \CY_B \: \vert \: \CY_{C \cup D}$. Similar problems arise for (P4) and (P5).
\end{itemize}
\end{remark}

The peculiarity of the partial correlation relation in continuous time is that it defines a graphoid under minimal assumptions. We require only wide-sense stationarity, zero expectation, mean-square continuity, and a positive definite spectral density function. For many graphoids, (P5) is quite difficult to verify and, unlike the proofs of (P1)--(P4),  additional, possibly strict, assumptions are required. For example, the conditional orthogonality relation for linear spaces  (\citealp{EI07}, Proposition A.1) satisfies (P5) only under the additional assumption of conditional linear separation of the underlying linear spaces. Thus, graphical models for stochastic processes using conditional orthogonality have additional assumptions on the spectral density (\citealp{EI07}, Eq. (2.1), for processes in discrete time and \citealp{VF23pre}, Assumption 1, for processes in continuous time) which guarantee that conditional linear separation holds.
Similarly, graphical models based on conditional independence also require additional assumptions (cf.~\citealp{LA04}, Proposition 3.1, and \citealp{EI11}, Assumption S).

\section{Partial correlation graphs}\label{sec: Partial correlation graphs}
First, in \Cref{subsec: Partial correlation graphs and Markov properties} we introduce the partial correlation graph $G_{PC}=(V, E_{PC})$, an undirected graph. This graph serves as a simple visual representation of the partial correlation structure within the multivariate stochastic process $\CY_V$. Moreover, for the partial correlation graph, we also derive edge characterisations and Markov properties. Finally, in \Cref{subsec: Partial correlation graphs and causality graphs}, we compare and contrast the partial correlation graph to the causality graph of \cite{VF23pre}. 

\subsection{Partial correlation graphs and Markov properties}\label{subsec: Partial correlation graphs and Markov properties}
Our approach to visualising the partial correlation structure between the components of the multivariate process $\CY_V$ in the graphical model $G_{PC}=(V, E_{PC})$ is as follows: Each component $\CY_a$, $a\in V$, is represented by a vertex. We then define a missing edge $a\edge b \notin E_{PC}$ if and only if the components $\CY_a$ and $\CY_b$ are uncorrelated given the linear information provided by  $\CY_{V \setminus \{a,b\}}$. As the relation $\CY_a \inde  \CY_b \: \vert \:  \CY_{V\setminus \{a,b\}}$ is symmetric, we use undirected edges in $G_{PC}$. This leads to the following definition of the partial correlation graph.

\begin{definition}\label{definition of the partial correlation graph}
Suppose that $\CY_V$ is wide-sense stationary with expectation zero, mean-square continuous, and has a spectral density function with $f_{Y_V Y_V}(\lambda)>0$ for $\lambda\in \R$. Let $V=\{1,\ldots,k\}$ be the vertices and define the edges $E_{PC}$ for $a,b\in V$ with $a\not=b$ as
\begin{align*}
    a \edge b \notin E_{PC}
    &\quad \Leftrightarrow \quad  \CY_a \inde  \CY_b \: \vert \:  \CY_{V\setminus \{a,b\}}.
\end{align*}
Then $G_{PC}=(V,E_{PC})$ is called \textsl{partial correlation graph} for $\CY_V$.
\end{definition}

\begin{remark}\mbox{}
\begin{itemize}
\item[(a)] The name partial correlation graph is clearly based on the partial correlation relation.
\item[(b)] For the definition of $G_{PC}$ it is not necessary to require that $f_{Y_V Y_V}(\lambda)>0$, but it is sufficient that $f_{Y_{V\setminus \{a,b\}}Y_{V\setminus \{a,b\}}}(\lambda)>0$ for all $a,b\in V$. However, $f_{Y_V Y_V}(\lambda)>0$ is essential for the graphoid properties and thus for the Markov properties of the partial correlation graph in \Cref{All markov properties satisfied}. Note that in general $f_{Y_V Y_V}(\lambda) \geq 0$ holds (cf.~\Cref{properties of spectral density}), so $f_{Y_V Y_V}(\lambda)>0$ is only a mild assumption.
\item[(c)] A direct consequence of \Cref{sufficient criterion for non-singular}(c) is that for $a\in V$ we would always have $a\edge a \in E_{PC}$. Since such self-loops do not help to visualise the partial correlation structure and do not change the properties of the graph, we omit them for the sake of simplicity.
\end{itemize}
\end{remark}

\begin{lemma}\label{characterisation of edges via density}
Suppose that $G_{PC} =(V, E_{PC})$ is the partial correlation graph for $\CY_V$. Then, for $a,b\in V$ with $a\neq b$, the following equivalences hold.
\begin{align*}
a \edge b \notin E_{PC}
&\quad \Leftrightarrow \quad  \cab(t) = 0 \quad \:\: \text{for all $t \in \R$,} \\
&\quad \Leftrightarrow \quad \fab(\lambda) = 0 \quad \: \text{for almost all $\lambda \in \R$,} \\
&\quad \Leftrightarrow \quad \Rab(\lambda) = 0 \quad \text{for almost all $\lambda \in \R$.}
\end{align*}
\end{lemma}

Note that the spectral coherence function is well-defined, since $f_{Y_VY_V}(\lambda)>0$ by assumption and thus, $\faa(\lambda) >0$ and $\fbb(\lambda) >0$, which results in a non-vanishing denominator.

In addition, \Cref{characterisation with inverse} gives the key representation using the inverse of the spectral density function, the corresponding edge characterisation for time series in discrete time is established in \cite{DA00}, Theorem 2.4.

\begin{proposition}\label{characterisation via inverse density}
Suppose that $G_{PC}=(V,E_{PC})$ is the partial correlation graph for $\CY_V$. Then, for $a,b\in V$ with $a\neq b$, the spectral coherence function satisfies
\begin{align*}
\Rab(\lambda) = - \frac{\left[g_{Y_VY_V}(\lambda)\right]_{ab}}{\left( \left[g_{Y_VY_V}(\lambda)\right]_{aa}\left[g_{Y_VY_V}(\lambda)\right]_{bb} \right)^{1/2}}, \quad \text{ $\lambda \in \R$.}
\end{align*}
Furthermore,
\begin{align*}
a \edge b \notin E_{PC} \quad \Leftrightarrow \quad \left[g_{Y_VY_V}(\lambda)\right]_{ab} = 0 \quad \text{for almost all $\lambda \in \R$. }
\end{align*}
\end{proposition}

\begin{remark}
A significant advantage of \Cref{characterisation via inverse density} over other characterisations is that it is computationally inexpensive. One only needs to know the spectral density function and then perform a singular matrix inversion to obtain all the edges in the graph simultaneously. 
Furthermore, the relation in \Cref{characterisation via inverse density} even gives us a simple measure for the strength of the connection between the components.
\end{remark}

\begin{remark}
\cite{DA00} remarks that the partial correlation graph can be compared to the concentration graph. The concentration graph for  a random vector $Z\in \R^k$ with $\BE \Vert Z \Vert^2<\infty$, $\BE(Z)=0_k$, and $ \Sigma_Z\coloneqq\BE [Z Z^\top ]>0$ is defined as follows. Let $V=\{1,\ldots,k\}$ be the vertices and define the edges $E_{CO}$ for $a,b\in V$ with $a\not=b$ as
\begin{align*}
    a\edge b \notin E_{CO} \quad &\Leftrightarrow \quad
\left[\BS_Z^{-1}\right]_{ab} =0\\
\quad &\Leftrightarrow \quad
    \left[\Sigma_Z\right]_{ab} - \left[\Sigma_Z\right]_{a V\setminus\{a,b\}}
    \left(\left[\Sigma_Z\right]_{V\setminus\{a,b\} V\setminus\{a,b\}}\right)^{-1}
    \left[\Sigma_Z\right]_{V\setminus\{a,b\} b}
    =0.
\end{align*}
 Then $G_{CO}=(V,E_{CO})$ is called \textsl{concentration graph}  of $Z$.
 The concentration graph $G_{CO}$ describes the sparsity pattern of the concentration matrix of $Z$. The definition of the concentration graph illustrates why the partial correlation graph for stochastic processes is a generalisation of the concentration graph for random vectors. A missing edge $ a\edge b \notin E_{CO}$  in the  concentration graph for $Z$ reflects that $Z_a$ and $Z_b$ are partially uncorrelated given  $Z_{V\setminus\{a,b\}}$. Similarly, $ a \edge b \notin E_{PC} $ in the partial correlation graph for $\CY_V$
 means that the stochastic processes $\CY_a$ and $\CY_b$ are partially uncorrelated given  $\CY_{V\setminus\{a,b\}}$. Finally, the edges in the partial correlation graph are characterised by the inverse of the spectral density function, which can be seen as a generalisation of the inverse of a covariance matrix. Indeed, for an independent and identically distributed sequence of random vectors with distribution $Z$, the spectral density is equal to $(2\pi)^{-1}\Sigma_Z$. 
    Note that the concentration graph is usually defined only for multivariate Gaussian random vectors \cite[p.~218]{Handbook:graphical} and not for general random vectors, but this definition is a natural generalisation. For Gaussian random vectors, however, missing edges correspond even to conditional independence relations \cite[Corollary 9.1.2]{Handbook:graphical}.
\end{remark}


To conclude this section, we establish the Markov properties of $G_{PC}$. To do this, we first provide some terminology.

\begin{definition}\label{def separation}
For $a\in V$ define $\ne(a)=\{ b \in V \: \vert \: a \edge b \in E_{PC} \}$ as the set of \textsl{neighbours} of $a$. A \textsl{path} of length $n$ from a vertex $a$ to a vertex $b$ is a sequence \linebreak $\alpha_0=a, \alpha_1,\ldots,\alpha_n=b$ of vertices such that $\alpha_{i-1} \edge \alpha_i \in E_{PC}$ for $i=1,\ldots,n$. For $A,B,C \subseteq V$, we say that $C$ \textsl{separates} $A$ and $B$ if every path from an element of $A$ to an element of $B$ contains at least one vertex from the separating set $C$. We write $A \bow B \: \vert \: C$ for short.
\end{definition}

Now the partial correlation graph satisfies the following Markov properties.

\begin{proposition}\label{All markov properties satisfied}
Suppose that $G_{PC}=(V,E_{PC})$ is the partial correlation graph for $\CY_V$. Then $\CY_V$ satisfies
\begin{itemize}
    \item[(P)] the \textsl{pairwise Markov property} with respect to $G_{PC}$, i.e., for $a,b\in V$ with $a\not= b$,
        \begin{itemize}
        	\item[]  $a\text{ \textbf{\---- }} b\notin E_{PC} \quad \Rightarrow \quad \CY_a \inde \CY_b \:\vert\: \CY_{V\setminus\{a,b\}}$, \phantom{$\Big[ \Big]$}
        \end{itemize}
    \item[(L)] the \textsl{local Markov property} with respect to $G_{PC}$, i.e., for $a\in V$,
        \begin{itemize}
			\item[] $\CY_{V\setminus (\ne(a)\cup \{a\})} \inde \CY_a \:\vert\: \CY_{\ne(a)}$, \phantom{$\Big[ \Big]$}
        \end{itemize}
     \item[(G)] the \textsl{global Markov property} with respect to $G_{PC}$, i.e., for disjoint $A,B,C \subseteq V$,
        \begin{itemize}
			\item[] $A \bow B \:\vert\: C \quad \Rightarrow\quad \CY_A \inde \CY_B \: \vert \: \CY_C$. \phantom{$\Big[ \Big]$}
        \end{itemize}
\end{itemize}
\end{proposition}

The pairwise Markov property holds by definition. Furthermore, the partial correlation relation defines a graphoid by \Cref{graphoid}. Thus, \cite{LA04} states in Theorem 3.7 that the pairwise, local and global Markov properties are equivalent, so the local and global Markov properties are also valid.
The global Markov property is important because it provides a graphical criterion for deciding when two subprocesses $\CY_A$ and $\CY_B$ are partially uncorrelated given a third subprocess $\CY_C$. Although the graph itself is defined only by pairwise partial correlation relations, we can obtain partial correlation relations between multivariate subprocesses given any subprocesses through path analysis.

\subsection{Partial correlation graphs and causality graphs}\label{subsec: Partial correlation graphs and causality graphs}
In this section, we draw parallels between the causality graph of \cite{VF23pre} and our partial correlation graph. First, we introduce the causality graph, using their edge characterisations in Lemmatas~3.2 and 4.2.

\begin{definition}\label{Definition causality graph}
Suppose that $\CY_V$ is wide-sense stationary with expectation zero, mean-square continuous, purely non-deterministic, and has a spectral density function with $f_{Y_V Y_V}(\lambda)>0$ for $\lambda\in \R$ that satisfies Assumption 1 of \cite{VF23pre}. Let $V=\{1,\ldots,k\}$ be the vertices and define the edges $E_{GC}$, for $a,b\in V$ with $a\neq b$, as
\begin{itemize}
\item[(i)\phantom{i}] \makebox[2,5cm][l]{$a \rarrow b \notin E_{GC}$}  
$\Leftrightarrow$ \quad $\CY_a$ is Granger non-causal for $\CY_b$ with respect to $\CY_V$
\item[\phantom{(ii)}]\makebox[2,5cm][l]{}  
$\Leftrightarrow$   $\quad {L}_{Y_b}(t+h) \perp \mathcal{L}_{Y_a}(t) \: \vert \: \mathcal{L}_{Y_{V\setminus\{a\}}}(t) \: \forall \: 0 \leq h \leq 1, \: t \in \R$,
\item[(ii)] \makebox[2,5cm][l]{$a \inst b \notin E_{GC}$}
$\Leftrightarrow$ \quad $\CY_a$ and $\CY_b$ are \textsl{contemporaneously uncorrelated} with respect
\item[\phantom{(ii)}]\makebox[3,5cm][l]{} to $\CY_V$
\item[\phantom{(ii)}]\makebox[2,5cm][l]{}   $\Leftrightarrow \quad L_{Y_a}(t+h) \perp L_{Y_b}(t+h') \: \vert \: \mathcal{L}_{Y_V}(t) \: \forall \: 0 \leq h, h' \leq 1, \: t \in \R$.
\end{itemize}
 Then $G_{GC}=(V,E_{GC})$ is called \textsl{(mixed) causality graph} for $\CY_V$. The index GC stands for \textsl{G}ranger \textsl{c}ausality.
\end{definition}

\begin{remark}\mbox{}
    To highlight the differences between the undirected edges in the causality graph and in the partial correlation graph, recall from \Cref{Bedingte Orthogonalitaet} that in the partial correlation graph, for $a,b\in V$ with $a \neq b$,
    \begin{align*}
    a \edge b \notin E_{PC}
    \quad \Leftrightarrow \quad
    L_{Y_a}(t) \perp L_{Y_b}(t) \: \vert \: \mathcal{L}_{Y_{V\setminus\{a,b\}}} \quad \forall \: t\in \R.
    \end{align*}
    The concept of contemporaneous uncorrelatedness in \Cref{Definition causality graph}(ii) differs from zero partial correlation in two ways. First, for zero partial correlation, we always project on the linear space of the whole process $\CY_{V\setminus\{a,b\}}=(Y_{V\setminus\{a,b\}}(t))_{t\in \R}$, whereas, for contemporaneous uncorrelatedness, we project on the past $(Y_V(s))_{s\leq t}$. Second, in the case of contemporaneous uncorrelatedness, the correlation has to be considered not only at identical time points but also at mixed time points one time step into the future.
\end{remark}

Despite the differences between the two concepts (which is also confirmed by the analysis of MCAR processes in  \Cref{Gegenbeispiel fuer Kantenbeziehungen}), there are relationships between the paths in the mixed causality graph and the edges in the partial correlation graph. To show these relations, we first provide the concept of $m$-separation (cf.~\citealp{EI07}), which is the extension of separation for undirected graphs (cf.~\Cref{def separation}) to mixed graphs.


\begin{definition}
In a mixed graph $G=(V, E)$ an intermediate vertex $c$ on a path $\pi$ is said to be a  \textsl{collider}, if the edges preceding and succeeding $c$ on the path both have an arrowhead or a dashed tail at $c$, i.e., $\rarrow c \larrow$, $\rarrow c \inst$, $\inst c \larrow $, or $\inst c \inst$. A path $\pi$ between vertices $a$ and $b$ is said to be \textsl{$m$-connecting} given a set $C$ if
\begin{itemize}
\item[(a)] every non-collider on $\pi$ is not in $C$, and
\item[(b)] every collider on $\pi$ is in $C$,
\end{itemize}
otherwise we say $\pi$ is \textsl{$m$-blocked} given $C$. If all paths between $a$ and $b$ are $m$-blocked given $C$, then $a$ and $b$ are said to be  \textsl{$m$-separated} given $C$, denoted by $ \{a\} \msep \{b\} \: \vert \: C  \:\: [G]$.
\end{definition}

The first relation between the causality graph and the partial correlation graph follows almost directly from the global AMP Markov property of the causality graph, which is established by \cite{VF23pre} in their Theorem 5.15.

\begin{lemma}\label{Zusammenhang zu causality graph}
Suppose that $G_{PC}=(V,E_{PC})$ is the partial correlation graph and $G_{GC}=(V,E_{GC})$ is the causality graph for $\CY_V$.  Then, for $a,b\in V$ with $a \neq b$, the following implication holds.
\begin{align*}
    \{a\} \msep \{b\} \: \vert\: V\setminus\{a,b\} \:\:\: [G_{GC}]
    \quad \Rightarrow \quad
    a\edge b \notin E_{PC}.
\end{align*}
\end{lemma}

The advantage of this result is that the concept of $m$-separation has several different characterisations in the literature, leading to more sufficient criteria for $a\edge b \notin E_{PC}$. One approach is to build an undirected graph from the mixed graph, using augmentation. The resulting augmented graph can then be related to the undirected partial correlation graph. The augmented graph is constructed as follows \cite[p.~148]{RI03}.


\begin{definition}
Let $G=(V, E)$ be a mixed graph. Two vertices $a$ and $b$ are said to be \textsl{collider connected} if they are connected by a pure collider path, which is a path on which every intermediate vertex is a collider. Then the undirected \textsl{augmented graph} $G^{a} = (V, E^{a})$ is derived from $G = (V, E)$ via
    \begin{align*}
        a\edge b \notin E^a
        \quad &\Leftrightarrow \quad
        \text{$a$ and $b$ are not collider connected in $G$.}
    \end{align*}
\end{definition}

Note that every single edge is trivially considered to be a collider path. Thus, every directed and undirected edge in the causality graph corresponds to an undirected edge in the augmented causality graph, implicating that the augmented causality graph has more edges than the causality graph. 

\begin{lemma}\label{Zusammenhang 2 zu causality graph}
Suppose that $G_{PC}=(V,E_{PC})$ is the partial correlation graph, \mbox{$G_{GC}=(V,E_{GC})$} is the causality graph, and $G_{GC}^a=(V,E_{GC}^a)$ is the augmented causality graph for $\CY_V$. For $a,b\in V$ with $a \neq b$, the following equivalences hold.
\begin{align}
    a\edge b \notin E^a_{GC}
    \quad &\Leftrightarrow \quad
    \dis\left( a \cup \ch(a)\right)\cap \dis\left( b \cup \ch(b) \right)=\emptyset \: \text{ in $G_{GC}$}, \label{Crit 1}\\
    \quad &\Leftrightarrow \quad
    \{a\} \bowtie \{b\} \: \vert\: V\setminus\{a,b\} \:\:\: [G_{GC}^a]. \label{Crit 2}
\end{align}
Here  $\ch(a)=\{v\in V \vert a \rarrow v \in E_{GC} \}$, $\dis(a)= \{ v\in V \vert v \inst \cdots \inst a \text{ or } v=a \}$ and
$\dis(A) = \bigcup_{a\in A} \dis(a)$. In particular, this implies that $a\edge b \notin E_{PC}$, i.e., $E_{PC} \subseteq E^{a}_{GC}$. 
\end{lemma}

This result gives us several possibilities to make statements about the partial correlation graph from the causality graph. On the one hand, the criterion \eqref{Crit 1} is particularly useful, since we can work with the original mixed graph and it is easy to implement algorithmically \citep{EI11}, which is not straightforward for the $m$-separation criterion from \Cref{Zusammenhang zu causality graph}. On the other hand, the characterisation \eqref{Crit 2} is of interest, as it uses the classical separation in an undirected graph, which is another common way to define global Markov properties in mixed graphs. Finally, the inclusion property $E_{PC} \subseteq E^{a}_{GC} $ gives us a simple connection between the edges in both graphs.



Besides the causality graph $G_{GC}=(V, E_{GC})$, \cite{VF23pre} also introduce the local causality graph $G_{GC}^0=(V, E_{GC}^0)$, a mixed graph with \mbox{$E_{GC}^0 \subseteq E_{GC}$}. For the augmented local causality graph obviously $E_{GC}^{a,0} \subseteq E_{GC}^a$ holds, but in general the statement $E_{PC} \subseteq E^{a,0}_{GC}$ is probably not possible, since we do not have a global AMP Markov property in the local causality graph.  However, if we restrict to MCAR$(p)$ processes, we derive this subset relation in \Cref{Partial correlation graphs and other graphical models}.

\section{Partial correlation graphs for MCAR processes}\label{sec: Partial correlation graphs for MCAR processes}
In the following, we construct the partial correlation graph for Lévy-driven multivariate continuous-time autoregressive (MCAR) processes to illustrate the partial correlation structure within this important and versatile class of processes. Therefore, in \Cref{subsec: MCAR processes}, we give a brief introduction to MCAR processes. Subsequently, in \Cref{subsec: The analysis of the assumptions for the MCAR process}, we ensure that the partial correlation graph for MCAR processes is well defined and establish the latter. We also provide some edge characterisations by model parameters along with comparisons to the literature. Moving on to \Cref{Partial correlation graphs and other graphical models}, we study relations between the partial correlation graph and the (local) causality graph, highlighting both similarities and differences.
Finally, in \Cref{sec:estimation}, we motivate some methods to estimate the edges in the partial correlation graph for MCAR processes.

\subsection{MCAR processes}\label{subsec: MCAR processes}
Early works on univariate and multivariate CAR processes and the more general continuous-time autoregressive moving average (CARMA) processes include those of \cite{doob:1944, DO60, Harvey:Stock:1985, Harvey:Stock:1988, HarveyStock1989, Bergstrom:1997}. Since then, these processes have enjoyed great popularity and have stimulated a considerable amount of research in recent years (cf.~\citealp{BR14}).
The driving process of an MCAR process is a $\R^k$-valued Lévy process $(L(t))_{t \in \R}$, which is a stochastic process with stationary and independent increments, it is continuous in probability, and satisfies \mbox{$L(0)=0_k \in \R^k$} $\mathbb{P}$-a.s.~ A typical example of Lévy process is the Brownian motion and the Poisson process. For more details on Lévy processes, we refer to the monographs of \cite{AP11} and \cite{ SA07}.
The following definition of a Lévy-driven MCAR process 
goes back to \cite{MA07}, Definition 3.20. 

\begin{definition}\label{Definition des MCAR Prozesses}
Let $L=(L(t))_{t\in \R}$ be a Lévy process satisfying $\BE[L(1)]=0_k$ and $\BE \Vert L(1) \Vert^2 < \infty$ with $\BS_L=\BE[L(1)L(1)^\top]$. Suppose that $A_1, A_2,\ldots, A_p \in \R^{k\times k}$ and define the matrices
    \begin{align*}
    \BA &= \begin{pmatrix}
      0_k & I_k & 0_k & \cdots & 0_k \\
      0_k & 0_k & I_k & \ddots & \vdots \\
      \vdots &  & \ddots & \ddots & 0_k \\
      0_k & \cdots & \cdots & 0_k & I_k \\
      -A_{p} & -A_{p-1} & \cdots & \cdots & -A_1
    \end{pmatrix}
    \in \R^{kp \times kp}, \quad
    \BB
    = \begin{pmatrix}
    0_k \\
    \vdots \\
    0_k \\
    I_k
    \end{pmatrix}
    \in \R^{kp\times k}, \\
    \BFC &= \begin{pmatrix}
       I_k & 0_k & \cdots & 0_k
    \end{pmatrix}
    \in \R^{k\times kp}.
    \end{align*}
Finally, suppose $\sigma(\BA)\subseteq (-\infty, 0) + i \R$ and $\CX=(X(t))_{t\in\R}$ is the unique $kp$-dimensional causal strictly stationary solution of the state equation
    \begin{align*}
    dX(t)= \BA X(t)dt+ \BB dL(t).
    \end{align*}
Then the output process $\CY_V=(Y_V(t))_{t\in\R}$ given by
\begin{align*}
 Y_V(t) = \BFC X(t)
\end{align*}
is called a \textsl{(causal) multivariate continuous-time autoregressive process of order $p$}, or MCAR$(p)$ process for short.
\end{definition}

The MCAR process is the continuous-time counterpart of the well-known discrete-time vector autoregressive (VAR) process.
For this correspondence, the idea is that a $k$-dimensional MCAR$(p)$ ($p\geq 1$) process $\CY_V$ is the solution to the stochastic differential equation
    \begin{equation*} \label{eq1.1}
     P(D) Y_V(t)=DL(t) 
    \end{equation*}
where $D$ is the differential operator with respect to $t$, and
    \begin{align} \label{ARpol}
     P(z) = I_k z^{p} + A_1 z^{{p}-1} + \ldots + A_{p}, \quad z\in \C,
    \end{align}
 is the autoregressive (AR) polynomial. However, a Lévy process is not differentiable, so this is not a formal definition of an MCAR process.
The properties of MCAR processes relevant to this paper are summarised below. For additional information, refer to \cite{MA07} and \cite{SC12, SC122}.

\newpage

\begin{remark}\label{MCAR wide-sense stationary}\mbox{}
\begin{itemize}
    \item[(a)] Since the input process $\CX$ is strictly stationary, the MCAR process $\CY_V$ is also strictly stationary. Furthermore, given the finite second moments of the Lévy process, both $\CX$ and $\CY_V$ also have finite second moments. 
    Thus, of course, the strictly stationary processes $\CX$ and $\CY_V$ are also wide-sense stationary.
    \item[(b)] The covariance function of the input process $\CX$ satisfies
    \begin{align*}
        c_{XX}(t)
        = \overline{c_{XX}(-t)}^\top
        = e^{\BA t} \Gamma(0), \quad t\geq 0,
        \quad \text{where} \quad
        \Gamma(0) = \int_{0}^{\infty} e^{\BA u} \BB \BS_L \BB^\top e^{\BA^\top u}du.
    \end{align*}
    The covariance function of the MCAR process $\CY_V$ is then determined as
    \begin{align*}
        c_{Y_V Y_V}(t) = \BFC c_{XX}(t) \BFC^\top,  \quad t \in \R.
    \end{align*}
    \item[(c)] Given that $\sigma(\BA)\subseteq (-\infty, 0) + i \R$, it follows that $c_{Y_VY_V}(t)$ decreases exponentially fast as \mbox{$t\rightarrow \pm \infty$} and
    $\lim_{t \rightarrow 0} c_{Y_VY_V}(t) = c_{Y_VY_V}(0)$, so $\CY_V$ is mean-square continuous due to \eqref{symmetrie of covariance}.
    \item[(d)] The MCAR$(1)$ process is also known as  Ornstein-Uhlenbeck process and in this case, we have $\BA=-A_1$ and $\BB = \BFC = I_k$. Furthermore, Gaussian MCAR processes and Gaussian Ornstein-Uhlenbeck processes, where the Brownian motion is the driving Lévy process, are special cases.
\end{itemize}
\end{remark}

\subsection{Definition of the partial correlation graph for a MCAR process} \label{subsec: The analysis of the assumptions for the MCAR process}
We introduce the partial correlation graph for MCAR processes. From  \Cref{MCAR wide-sense stationary} we already know that the MCAR process is wide-sense stationary with expectation zero and mean-square continuous. Furthermore, the spectral density function is \cite[Eq.~(3.43)]{MA07}
    \begin{align*}
        f_{Y_VY_V}(\lambda)
        &=\frac{1}{2\pi} P(i\lambda)^{-1} \BS_L \left (P(-i\lambda)^{-1} \right)^\top,
        \quad \lambda \in \R,
    \end{align*}
where the AR polynomial $P$ is defined in equation \eqref{ARpol}. For the well-definedness of the partial correlation graph we then only need to ensure that $f_{Y_V Y_V}(\lambda)>0$ for $\lambda \in \R$. But this condition is already met when $\BS_L > 0$ and $\sigma(\BA)\subseteq (-\infty, 0) + i\R$. Then the  inverse spectral density function has the representation
\begin{align*}
g_{Y_V Y_V}(\lambda)
= 2\pi P(-i\lambda)^\top \BS_L^{-1} P(i\lambda), \quad \lambda \in \R.
\end{align*}
By \Cref{definition of the partial correlation graph}, \Cref{characterisation via inverse density}, and \Cref{All markov properties satisfied} we then obtain the following result.

\begin{proposition}\label{Prop: Definition partial correlation graph for MCAR}
Suppose $\CY_V$ is a causal MCAR$(p)$ process with \mbox{$\BS_L>0$.} Let \linebreak $V=\{1,\ldots,k\}$ be the vertices and define the edges $E_{PC}$, for $a,b\in V$ with $a\not= b$, via
\begin{align*}
a \edge b \notin E_{PC}
\quad \Leftrightarrow \quad
\CY_a \inde \CY_b \: \vert \: \CY_{V\setminus \{a,b\}}
\quad \Leftrightarrow \quad
\left[ P(-i\lambda)^\top \BS_L^{-1} P(i\lambda) \right]_{ab}=0 \quad  \forall \: \lambda\in \R.
\end{align*}
 Then the partial correlation graph $G_{PC}=(V, E_{PC})$ for the MCAR process $\CY_V$ is well defined and satisfies the pairwise, local, and global Markov property.
\end{proposition}

Note that partial correlation graphs can be defined for more general state space models, but we find that MCAR processes are sufficient for our illustrative purposes. Note also that for the MCAR process, $g_{Y_V Y_V}(\lambda)$ has a very simple representation, it is a matrix polynomial. As a result, we can give the following edge characterisation based on the coefficients of the matrices $A_1, A_2, \ldots, A_p$ of the AR polynomial, and the covariance matrix $\Sigma_L$ of the driving Lévy process.

\begin{proposition}\label{Characterisierung PC fuer MCAR}
Suppose that $G_{PC}=(V, E_{PC})$ is the partial correlation graph for the causal MCAR$(p)$ process $\CY_V$ with AR polynomial $P$ given by \eqref{ARpol}, where we define $A_0:=I_k$. For $a,b\in V$ with $a\not= b$, we obtain the edge characterisation 
\begin{align*}
a\edge b \notin E_{PC}
\quad \Leftrightarrow \quad
 \sum_{m=0\vee n-p}^{n\wedge p} (-1)^m \left[ A^\top _{p-m}\BS_L^{-1}A_{p-n+m}
\right]_{ab}=0 \quad \text{for $n=0,\ldots,2p$.}
\end{align*}
This characterisation is reduced in the following cases.
\begin{itemize}
    \item[(i)] Suppose $\BS_L=\sigma^2 I_k>0$. Then
\begin{align*}
a\edge b \notin E_{PC}
\quad \Leftrightarrow \quad
 \sum_{m=0\vee n-p}^{n\wedge p} (-1)^m \left[A^\top _{p-m}A_{p-n+m}
\right]_{ab}=0 \quad  \text{for $n=0,\ldots,2p$.}
\end{align*}
    \item[(ii)] Suppose $A_j$ is a diagonal matrix for $j=1,\ldots,p$. Then
\begin{align*}
a \edge b \notin E_{PC}
\quad \Leftrightarrow \quad
\left[\BS_L^{-1}\right]_{ab} =0.
\end{align*}
\end{itemize}
\end{proposition}

\begin{remark}
 A consequence of  \Cref{Characterisierung PC fuer MCAR}(ii) is that for any undirected graph $G=(V, E)$ and any $p\in\N$, there exists an MCAR$(p)$ process with partial correlation graph $G_{PC}=G$. 
 Indeed, we can define
\begin{align*}
    \left[ \BS_L^{-1} \right]_{ab}
    = \begin{cases}
        k, & \text{if } a=b,\\
        1, & \text{if } a \neq b \text{ and } a\edge b \in E,\\
        0, & \text{if } a \neq b \text{ and } a\edge b \notin E,
    \end{cases}
\end{align*}
and 
$A_m=\left(p\atop m\right)I_k\in \R^{k\times k}$ for $m=0,\ldots,p$. Consequently, $\sigma(\BA)=\{-1\} \subseteq (-\infty, 0)+i\R$ and $\BS_L^{-1}$ is strictly diagonally dominant, i.e., positive definite. $\BS_L$ is also positive definite and there exists a Lévy process with this covariance matrix. Due to   \Cref{Characterisierung PC fuer MCAR}(ii) the resulting $k$-dimensional MCAR$(p)$ process $\CY_V$ generates a partial correlation graph $G_{PC}=(V,E_{PC})$, which is identical to the undirected graph $G=(V,E)$. This is a major advantage over the causality graph in \cite{VF23pre}, where it is not clear if any graph can be constructed by a continuous-time process.
\end{remark}

\begin{remark} \label{Remark 5.6}
The edge characterisations for MCAR$(p)$ processes in \Cref{Characterisierung PC fuer MCAR} are, as might be expected, similar to the edge characterisations for VAR$(p)$ processes in \cite{DA00}, Example 2.2. Suppose that the AR coefficient matrices of the VAR$(p)$ process are denoted by $\Phi_m\in \R^{k \times k}$, $m=1,\ldots,p$, $\Phi_0=-I_k$, and  $0<\Sigma_{\varepsilon}\in \R^{k \times k}$ denotes the covariance matrix of the white noise process. Then \cite{DA00} states that in the  partial correlation graph $G_{PC}^d=(V,E_{PC}^d)$ for the VAR$(p)$ process we have
\begin{align*}
    a\edge b \notin E_{PC}^d
    \quad &\Leftrightarrow \quad
    \sum_{m=0 \vee n-p}^{p \wedge n} \left[\Phi_{m}^\top \Sigma_{\varepsilon}^{-1} \Phi_{m-n+p}\right]_{ab}= 0 \quad \text{for } n=0,\ldots,2p.
\end{align*}
Both characterisations of the continuous-time and the discrete-time multivariate AR processes match exactly if we neglect the factor $(-1)^m$. This small difference is due to the fact that the spectral density of the continuous-time model is defined by the AR polynomial at $\pm i\lambda$ whereas, in the discrete-time model, it is the AR polynomial at $e^{\pm i\lambda}$.
%
\end{remark}

Furthermore, the following sufficient condition for an edge between $a$ and $b$ in the partial correlation graph can be obtained by setting $n=2p$ in \Cref{Characterisierung PC fuer MCAR}.

\begin{lemma}\label{Sufficient fuer Kante}
Suppose that $G_{PC}=(V,E_{PC})$ is the partial correlation graph for the causal MCAR$(p)$ process $\CY_V$. For $a,b\in V$ with $a\not= b$, the following implication holds.
\begin{align*}
a\edge b \notin E_{PC}
\quad \Rightarrow \quad
\left[\BS_L^{-1}\right]_{ab} =0.
\end{align*}
\end{lemma}

\begin{remark}
Note that $\Sigma_L^{-1}$ is the concentration matrix of the random vector $L(1)$, so it defines the concentration graph $G_{CO}=(V, E_{CO})$ of $L(1)$. \Cref{Sufficient fuer Kante} therefore gives the subset relation $E_{CO}\subseteq E_{PC}$. In other words, the partial correlation of the random variables $L_a(1)$ and $L_b(1)$  given  $L_{V\setminus\{a,b\}}(1)$ imply an edge in the partial correlation graph of the continuous-time process $\CY_V$, i.e., the stochastic processes $\CY_a$ and $\CY_b$  are partially correlated given the process $\CY_{V\setminus\{a,b\}}$. If we additionally assume that $A_m$, \linebreak $m=1,\ldots,p$, are diagonal, then \Cref{Characterisierung PC fuer MCAR}(ii) even gives $E_{CO} = E_{PC}$.
\end{remark}

Finally, for a visualisation of the previous edge characterisations in \Cref{Prop: Definition partial correlation graph for MCAR} and \Cref{Characterisierung PC fuer MCAR}, we present an example.

\begin{example}\label{Example Ornstein Uhlenbeck Teil 1}
Suppose that $\CY_V$ is a 4-dimensional Ornstein-Uhlenbeck process with $\BS_L=I_4$ and
\begin{align*}
\BA=\left(\begin{array}{rrrr}
-2 & 0  & 1  &  1 \\
0  & -2 & -1 & -1 \\
-1 & -1 & -2 & -1 \\
1  & -1 & -1 & -2
\end{array}\right).
\end{align*}
Then a simple calculation yields $\sigma(\BA) = \{-1,-1,-2,-4\} \subseteq (-\infty, 0) + i\R$. For an Ornstein-Uhlenbeck process $\CY_V$ the inverse spectral density function is simplified to
$g_{Y_V Y_V}(\lambda) = 2 \pi (-i \lambda I_k-\BA^\top) \BS_L^{-1} (i \lambda I_k-\BA)$ for $\lambda \in \R$ and we obtain
\begin{gather*}
g_{Y_VY_V}(\lambda)
=2\pi \left(\begin{matrix}
\lambda^2+6      & 0           & 2 i \lambda -1 & -3 \\
0                & \lambda^2+6 & 5              & 5 \\
- 2 i \lambda -1 & 5           &  \lambda^2+7   & 6 \\
-3               & 5           & 6              &\lambda^2+7
\end{matrix}
\right).
\end{gather*}
The corresponding partial correlation graph $G_{PC}=(V,E_{PC})$ is then given in \Cref{Partial correlation graph for Example Proposition characterisation for MCAR edges 2}.
\begin{figure}[ht]
\begin {center}
\begin{tikzpicture}[align=center, node distance=1cm and 2cm, semithick, state/.style ={circle, draw,  text=black, minimum width =0.5 cm}, every loop/.style={}]

  \node[state] (1) {1};
  \node[state] (2) [right of=1] {2};
  \node[state] (3) [below of=1] {3};
  \node[state] (4) [right of=3] {4};

  \path
   (1) edge [] node {} (4)
   (2) edge [bend left] node {} (4)
   (1) edge [bend right] node {} (3)
   (2) edge [] node {} (3)
   (3) edge [bend right] node {} (4);

\end{tikzpicture}
\end{center}
\caption{Partial correlation graph for \Cref{Example Ornstein Uhlenbeck Teil 1}}
\label{Partial correlation graph for Example Proposition characterisation for MCAR edges 2}
\end{figure}
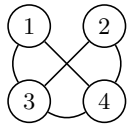

Furthermore, for an Ornstein-Uhlenbeck process, the edge characterisation in
\Cref{Characterisierung PC fuer MCAR}(i) is simplified to
\begin{align} 
  a \edge b \notin E_{PC}
\quad &\Leftrightarrow \quad
\left[\BA \right]_{ba} - \left[\BA \right]_{ab}=0,  \:\:\:
\left[ \BA^\top \BA \right]_{ab} =0. \label{OU Kriterium 1b}
\end{align}
Of course, this relation 
also provides the edges in \Cref{Partial correlation graph for Example Proposition characterisation for MCAR edges 2}. 
\end{example}

To summarise, \Cref{Example Ornstein Uhlenbeck Teil 1} highlights once more the main advantage of the characterisation in \Cref{Prop: Definition partial correlation graph for MCAR}, which is the ability to obtain all edges simultaneously through the inverse spectral density function.

\subsection{Partial correlation graphs and (local) causality graphs}\label{Partial correlation graphs and other graphical models}
In this section, we relate the partial correlation graph to the causality graph and the local causality graph of \cite{VF23pre}, which can be seen as a continuation of \Cref{subsec: Partial correlation graphs and causality graphs}. Let us start with the relations between the partial correlation graph and the \textsl{causality graph}. In the comparison in \Cref{subsec: Partial correlation graphs and causality graphs}, we suspected that, in general, there is no direct relationship between the edges in the causality graph and the partial correlation graph, although $E_{PC} \subseteq E_{GC}^a $. We now confirm this conjecture with two counterexamples.

\begin{example}\mbox{}\label{Gegenbeispiel fuer Kantenbeziehungen}
Recall that for the Ornstein-Uhlenbeck process with $\Sigma_L=I_k$, due to \Cref{Characterisierung PC fuer MCAR} with $p=1$, the characterisation
\begin{align*}
  a \edge b \notin E_{PC}
\quad &\Leftrightarrow \quad
\left[ \BA\right]_{ba} - \left[\BA \right]_{ab} =0, \quad
\left[ \BA^\top \BA \right]_{ab} =0,
\end{align*}
holds. Additionally, by Corollary 6.21 of \cite{VF23pre}, we have
\begin{align} \label{eq5.3}
\begin{array}{lll}
a\rarrow b \notin E_{GC}
\quad &\quad \Leftrightarrow \quad
\left[\BA^\alpha\right]_{ba}=0,  & \quad \alpha=1,\ldots,k-1, \\
a \inst b \notin E_{GC}
\quad &\quad\Leftrightarrow \quad
\left[\BA^\alpha \left(\BA^\top \right)^\beta \right]_{ab}=0,  & \quad \alpha,\beta=0,\ldots,k-1.
\end{array}
\end{align}

\begin{itemize}
    \item[(a)] Suppose that $\CY_V$ is a 3-dimensional Ornstein-Uhlenbeck process with $\Sigma_L=I_3$ and
\begin{align*}
    \BA=\left(\begin{array}{rrr}
        -3 & 1 & 1 \\
        1 & -3 & 1 \\
        6 & 1 & -8
    \end{array}\right),
\end{align*}
where $\sigma(\BA) = \{-9,-4,-1\} \subseteq (-\infty, 0) + i\R$. Then
\begin{align*}
\left[\BA\right]_{21} - \left[\BA \right]_{12} = 0 \quad \text{and} \quad
\left[ \BA^\top \BA \right]_{12}  = 0,
\end{align*}
so $1 \edge 2 \notin E_{PC}$. However $1 \rarrow 2 \in E_{GC}$, $2 \rarrow 1 \in E_{GC}$, and $1 \inst 2 \in E_{GC}$, since
\begin{align*}
   \left[\BA \right]_{21}  \neq 0  \quad \text{and} \quad
   \left[\BA \right]_{12}  \neq 0.
\end{align*}
    \item[(b)] Suppose that $\CY_V$ is a 3-dimensional Ornstein-Uhlenbeck process with \mbox{$\Sigma_L=I_3$ and}
\begin{align*}
    \BA=\left(\begin{array}{rrr}
        -1 & 0 & 0 \\
        0 & -1 & 0 \\
        1 & 1 & -2
    \end{array}\right),
\end{align*}
where $\sigma(\BA) = \{-1,-1,-2\} \subseteq (-\infty, 0) + i\R$. Then a simple calculation shows that
    \begin{align*}
        \left[ \BA^\alpha \right]_{2 1}= \left[ \BA^\alpha \right]_{1 2} =0, \quad \alpha=1,2, \quad \text{and} \quad
        \left[ \BA^\alpha \left(\BA^\top \right)^\beta \right]_{12}=0, \quad  \alpha,\beta=0,1,2.
    \end{align*}
 Therefore, $1 \rarrow 2 \notin E_{GC}$, $2 \rarrow 1 \notin E_{GC}$, and $1 \inst 2 \notin E_{GC}$. However, $1 \edge 2 \in E_{PC}$, since
   $        \left[ \BA^\top \BA \right]_{12}=1.$
\end{itemize}
\end{example}

In summary, even in the special case $\Sigma_L=I_k$, there are no direct relations between the edges because, in the partial correlation graph the orthogonality of the columns in $\BA$ is characteristic, whereas in the causality graph the orthogonality of the rows is relevant for the undirected edges, and the orthogonality of the rows and columns is relevant for the directed edges. Of course, in some special cases, there are simple relations between the edges in the partial correlation graph and the edges in the causality graph. Because of the orthogonality argument, an obvious special case is a symmetric matrix $\BA$.

\begin{lemma}\label{Spezialfall Beziehung}
Suppose that $G_{PC}=(V,E_{PC})$ is the partial correlation graph and $G_{GC}=(V,E_{GC})$ is the causality graph for the causal Ornstein-Uhlenbeck process $\CY_V$, where $\BA$ is a symmetric matrix and $\Sigma_L=I_k$.  Then, for $a,b\in V$ with $a\not=b$, we receive
\begin{align*}
    a \inst b \notin E_{GC} \quad \Rightarrow \quad a \edge b \notin E_{PC}.
\end{align*}
\end{lemma}

Next, we provide a comparison to the \textsl{local causality graph} established by \cite{VF23pre}. To avoid going too deep into the intricate definition of the local causality graph in its generality here, we present the definition of the local causality graph only for MCAR processes via the characterisations used in \cite{VF23pre}, Propositions 6.12 and 6.13. For a general definition of the local causality graph, we refer to their Definition 5.9.

\begin{definition}
Suppose $\CY_V$ is a causal MCAR$(p)$ process with $\BS_L>0$. Suppose $V=\{1,\ldots,k\}$ are the vertices and the edges $E_{GC}^0$ for $a,b\in V$ with $a\neq b$ are defined via
\begin{itemize}
\item[(i)] \
$a\rarrow b \notin E_{GC}^0
\quad \Leftrightarrow \quad
\left[A_j\right]_{ba}=0 \quad \text{for } j=1,\ldots,p,$
\item[(ii)] $
a \inst b \notin E_{GC}^0
\quad \Leftrightarrow \quad
\left[\Sigma_L\right]_{ab}=0.$
\end{itemize}
Then $G_{GC}^0=(V,E_{GC}^0)$ is called \textsl{local causality graph} for $\CY_V$.
\end{definition}

\begin{remark}\mbox{}
\begin{itemize}
\item[(a)] We emphasise that the undirected edges in the local causality graph are characterised by $\Sigma_L$ and not $\Sigma_L^{-1}$ as in the partial correlation graph, and these matrices generally do not match. The local causality graph considers the direct correlation of $L_a(1) $ and $L_b(1) $, while the partial correlation graph considers the correlation of $L_a(1) $ and $L_b(1) $ given the environment $L_{V\setminus \{a,b\}}(1)$.
\item[(b)]  Due to the different definitions, there are generally no direct relations between the edges in the partial correlation graph and the edges in the local causality graph, not even in the special case $\Sigma_L=I_k$. Note that in this case, $a\inst b \notin E_{GC}^0$ is always true. Furthermore, looking at \Cref{Gegenbeispiel fuer Kantenbeziehungen}(a), we get $1 \edge 2 \notin E_{PC}$ but $1 \rarrow 2 \in E_{GC}^0$ and $2 \rarrow 1 \in E_{GC}^0$. Whereas \Cref{Gegenbeispiel fuer Kantenbeziehungen}(b) is an example where $1 \rarrow 2 \notin E_{GC}^0$ and $2 \rarrow 1 \notin E_{GC}^0$ but $1 \edge 2 \in E_{PC}$.
\item[(c)] In the case of no environment ($k=2$) we obtain that $\left[\Sigma_L\right]_{ab}=0$ if and only if $[\Sigma_L^{-1}]_{ab}=0$ and $a\edge b \notin E_{PC}$ implies $a \inst b \notin E_{GC}^0$ and vice versa.  \qedhere
\end{itemize}
\end{remark}

However, as for the causality graph, we can establish relations between edges in the partial
correlation graph and paths in the local causality graph for MCAR processes via the concept of $m$-separation and augmentation separation, although no global AMP Markov property could be shown for the local
causality graph.

\begin{lemma}\label{Separationsvergleich local}
Suppose that $G_{PC}=(V,E_{PC})$ is the partial correlation graph, \mbox{$G_{GC}^0=(V,E_{GC}^0)$} is the local causality graph, and
$G_{GC}^{0,a}=(V,E_{GC}^{0,a})$ is the augmented local causality graph for the causal MCAR$(p)$ process $\CY_V$.  Then,
for $a,b\in V$ with $a \neq b$, the following equivalences hold.
\begin{align*}
 a\edge b \notin E^{0,a}_{GC}
     \quad &\Leftrightarrow \quad
     \{a\} \bowtie \{b\} \: \vert\: V\setminus\{a,b\} \:\:\: [G_{GC}^{0,a}],\\
     \quad &\Leftrightarrow \quad
     \{a\} \msep \{b\} \: \vert\: V\setminus\{a,b\} \:\:\: [G_{GC}^{0}], \\
     \quad &\Leftrightarrow \quad
     \dis\left( a \cup \ch(a)\right)\cap \dis\left( b \cup \ch(b) \right) \: \text{in $G_{GC}^0$}.
\end{align*}
In particular, we then have $a\edge b \notin E_{PC}$, i.e., $E_{PC} \subseteq E_{GC}^{0,a}$. 
\end{lemma}

Note that the opposite inclusion does in general not hold, there exist examples where $E_{PC}\not= E_{GC}^{0,a}$ as for the causality graph.

As discussed in \Cref{Zusammenhang 2 zu causality graph}, \Cref{Separationsvergleich local} provides us with several ways to make statements about the partial correlation graph from the local causality graph.

\subsection{Estimation} \label{sec:estimation}

 The edges in partial correlation graph can be found simultaneously and computationally inexpensive using the inverse of the spectral density function (cf.~\Cref{characterisation via inverse density}). Therefore, in practical applications, we have to estimate the spectral density function from discrete-time observations.
 Suppose we  observe a causal  MCAR$(p)$ process $\CY_V$ at equidistant times $0,\,\Delta,\,2\Delta,\ldots$ with $\Delta>0$ small, as used for modelling high-frequency data. The resulting discrete-time process $\CY_V^{\Delta}=(Y_V(k\Delta))_{k\in\N}$ is also weakly stationary with zero expectation, in fact, it is  a vector ARMA process, with spectral density function
\begin{equation*} \label{Eq 2}
 f_{Y_VY_V}^{(\Delta)}(\lambda)=\frac{1}{2\pi}\sum_{k=-\infty}^\infty c_{Y_VY_V}(k\Delta)\,{\rm e}^{-ik\lambda} = \frac{1}{\Delta}\sum_{k=-\infty}^\infty f_{Y_VY_V}\left(\frac{\lambda+2k\pi}{\Delta}\right),\quad-\pi\leq\lambda\leq\pi,
\end{equation*}
where the second equality follows from \cite{Bloomfield1976}, p.~206.

\subsubsection*{Low-frequency sampling scheme}

But clearly the zero entries of the inverse of $f_{Y_VY_V}^{(\Delta)}(\lambda)$ for $\lambda\in[-\pi,\pi]$
do not necessarily coincide with the zero entries of the inverse of $f_{Y_VY_V}(\lambda)$ for $\lambda \in\R$ and hence, there is in general no relationship between the edges in the partial correlation graph for $\CY_V$ and the partial correlation graph for $\CY_V^{(\Delta)}$. This can be seen nicely by looking at an Ornstein-Uhlenbeck process, where, due to \Cref{Characterisierung PC fuer MCAR}, we have
\begin{align*}
  a \edge b \notin E_{PC}
\quad &\Leftrightarrow \quad
\left[\Sigma_L^{-1} \right]_{ab}=0, \quad
\left[ \BA^\top \Sigma_L^{-1} - \Sigma_L^{-1} \BA \right]_{ab} =0, \quad
\left[ \BA^\top \Sigma_L^{-1} \BA \right]_{ab} =0.
\end{align*}
The discrete-time sampled process $\CY_V^{(\Delta)}$ is a VAR(1) process  where, due to \Cref{Remark 5.6},
the edges in the partial correlation graph  $G^{(\Delta)}_{PC}=(V,E^{(\Delta)}_{PC})$ for $\CY_V^{(\Delta)}$
can be described by the relation
\begin{align*}
a \edge b \notin E^{(\Delta)}_{PC}
\quad \Leftrightarrow \quad
&\left[ \left( \BS^{(\Delta)} \right)^{-1}  + e^{\BA ^\top \Delta} \left( \BS^{(\Delta)} \right)^{-1}  e^{\BA \Delta} \right]_{ab}=0, \\
&\left[ \left( \BS^{(\Delta)} \right)^{-1} e^{\BA \Delta} \right]_{ab}=0, \quad
\left[ e^{\BA ^\top \Delta} \left( \BS^{(\Delta)} \right)^{-1} \right]_{ab}=0
\end{align*}
where $\BS^{(\Delta)}=\int_0^\Delta e^{\BA u} \BS_L e^{\BA^\top u} du>0$.
These characterisations confirm that there do not exist direct relationships between $E_{PC}$ and $E^{(\Delta)}_{PC}$. Therefore, in general, it will be challenging to derive a nonparametric estimator for the edges in the partial correlation graph from a  low-frequency sampling scheme.

However, it is possible to derive estimators for a parametric class of continuous-time processes.
In the case of MCAR models, the model parameters $A_1,\ldots,A_p,\Sigma_L$ can be estimated from the low-frequency sampled MCAR process $\CY_V^{(\Delta)}$, e.g.,
by quasi maximum-likelihood estimation as in \cite{SC12} or Whittle estimation as in \cite{FHM}, yielding the parameter estimators $\widehat A_1,\ldots,\widehat A_p,\widehat \Sigma_L$, which are consistent and asymptotically normally distributed. Then an estimator for the inverse of the spectral density is
\begin{align*}
    \widehat g_{Y_V Y_V}(\lambda)&= 2\pi \widehat P(-i\lambda)^\top \widehat \BS_L^{-1} \widehat P(i\lambda), \quad \lambda \in \R \quad \text{with}\\
     \widehat P(z) &= I_k z^{p} + \widehat A_1 z^{{p}-1} + \ldots + \widehat A_{p}, \quad z\in\C,
\end{align*}
which is also a consistent and asymptotically normally distributed estimator for $g_{Y_V Y_V}(\lambda)$ for fixed $\lambda\in\R$ by an application of the continuous mapping theorem and the delta-method, respectively.
By considering the zero entries of this function we receive estimators for the edges in the partial correlation graph for the underlying continuous-time process $\CY_V$. 

\subsubsection*{High-frequency sampling scheme}

In the context of high-frequency data  where \mbox{$\Delta\to 0$}, we have the relation
\begin{equation} \label{Equation f}
 \lim_{\Delta\to 0}\,\Delta\,f^{\Delta}_{Y_VY_V}(\lambda\Delta)\,\mathds{1}_{[-\frac{\pi}{\Delta},\frac{\pi}{\Delta}]}(\lambda)= f_{Y_VY_V}(\lambda), \quad \lambda\in\R,
\end{equation}
(\citealp{FAFU13}, Eq. (1.5) for CARMA processes but this is also true for our causal MCAR processes). Roughly speaking, this means that in the limit $\Delta\to 0$, we can identify the edges of the causal MCAR process from edges of its equidistantly sampled observations.
In the special case of univariate CARMA processes, we already know from \cite{FAFU13, FAFU13b} that, under some mild assumptions, the  smoothed normalised periodogram  is a consistent estimator of the spectral density $f_{Y_VY_V}(\lambda)$ for the high-frequency sampling scheme, where $\Delta_n\to 0$
and $n\Delta_n\to\infty$ as the number of observations $n\to\infty$. We believe it is straightforward to show that this is still true for multivariate CARMA processes including MCAR processes. An alternative estimator is the lag-window spectral density estimator of \cite{kartsioukas2023spectral}. They develop the statistical inference of this estimator not only for MCAR processes but also for general multivariate stationary processes in Hilbert spaces and, furthermore, they also allow an irregular sampling scheme. For non-Gaussian processes, however, a cumulant condition must be satisfied, which is in the context of MCAR processes a cumulant condition on the driving Lévy process. Due to the generality of this impressive paper,  the assumptions for MCAR processes are actually stronger than necessary.

\section{Conclusion} \label{Conclusion}
The paper establishes and analyses the partial correlation relation for wide-sense stationary and mean-square continuous stochastic processes in continuous time with expectation zero and spectral density function. Based on this, the partial correlation graph for continuous-time stochastic processes is defined, which satisfies the usual Markov properties. Furthermore, we relate the partial correlation graph to the causality and the local causality graph by \cite{VF23pre} by their augmented graphs and we find some interesting relationships. The derived results for the partial correlation graph in the continuous-time setting correspond to the results for discrete-time processes in \cite{DA00}, which we also see by applications to MCAR processes, where we can characterise the edges by the model parameters. In both settings, the low-frequency sampling regime and the high-frequency sampling regime, it is possible to derive some consistent and asymptotically normally distributed estimators for the inverse spectral density of an MCAR process and thus also for the edges in the partial correlation graph. In the high-frequency sampling scheme, the smoothed periodogram and the lag-window spectral density estimator are popular estimators for the spectral density as for discrete-time processes (\citealp{Anderson71, BR91, BR01, Hannan1970}) and they should also work for a large class of non-parametric continuous-time models. The paper focused on the theoretical properties of the partial correlation graph but
 statistical methods for estimation and testing for the edges in the continuous-time partial correlation graph are of particular importance and will be the subject of some future work.

\bibliographystyle{imsart-nameyear}
\bibliography{102_Literatur}

\appendix
\section{Proofs}\label{sec:proofs}

\subsection{Proofs of Section~\ref{sec: Partial correlation}}\label{sec: Proofs of sec: Partial correlation}
\begin{proof}[Proof of \Cref{Gleichheit linearer Raeume}]
The relation $\subseteq$ in \Cref{Gleichheit linearer Raeume} is obvious, since we have $Y_c(t)=\uint e^{i \lambda t} \Phi_c(d\lambda)\in \mathcal{L}_{Y_C}^*$ for all $c\in C$ and $t\in \R$, and $\mathcal{L}_{Y_C}^*$ is a closed linear space. The relation $\supseteq$ is established by \cite{RO67} on p.~34.
\end{proof}

\begin{proof}[Proof of \Cref{Loesung ist Projektion}]
Let $t\in \R$ and assume that $\{a\}\cap C = \emptyset$, since the statements apply trivially for $a \in C$. To simplify the notation, we abbreviate
\begin{align*}
  \widehat{Y}_{a\vert C}(t) = \uint e^{i \lambda t} f_{Y_aY_C}(\lambda) f_{Y_CY_C}(\lambda)^{-1} \Phi_C(d\lambda).
\end{align*}
The proof is divided into three steps. In the first step we derive that $\widehat{Y}_{a\vert C}(t) \in \mathcal{L}_{Y_C}$ and in the second step we show that  $Y_a(t) - \widehat{Y}_{a\vert C}(t) \in \mathcal{L}_{Y_C}^\perp$. Both together then give the assertion $\widehat{Y}_{a\vert C}(t)=P_{\mathcal{L}_C}Y_a(t)$. Then, in a third step, we conclude that the orthogonal projection is the solution to the optimisation problem \eqref{optimisation problem}.\\
\textsl{Step 1: \:} Given that $\mathcal{L}_{Y_C}=\mathcal{L}_{Y_C}^*$ due to \Cref{Gleichheit linearer Raeume}, we can establish the measurability and integrability of the function $f_{Y_a Y_C}(\lambda) f_{Y_CY_C}(\lambda)^{-1}$ for $\lambda \in \R$.
For the measurability, we first note that $f_{Y_aY_C}$ and $f_{Y_CY_C}$ are 
measurable as derivatives. Furthermore, sums and products of measurable functions are measurable. If we set $\lambda/0\coloneqq0$ for $\lambda \in \R$, then their quotients are also measurable \cite[Theorem 1.91]{KL20}. Now we compute $f_{Y_CY_C}(\lambda)^{-1}$ by Gaussian elimination and find that $f_{Y_CY_C}(\lambda)^{-1}$ is measurable for $\lambda \in \R$. Thus, $f_{Y_a Y_C}(\lambda)f_{Y_CY_C}(\lambda)^{-1}$, $\lambda \in \R$, is also measurable. \\
For the integrability, we first note that $f_{Y_{\{a\} \cup C} Y_{\{a\} \cup C}} (\lambda) \geq 0$ due to \Cref{properties of spectral density}(c). Furthermore, $f_{Y_{C} Y_{C}} (\lambda) > 0$ by assumption, so Proposition 8.2.4 of \cite{BE09} gives
\begin{align*} 
    f_{Y_aY_C}(\lambda) f_{Y_CY_C}(\lambda)^{-1} f_{Y_CY_a}(\lambda) \leq f_{Y_aY_a}(\lambda)
\end{align*}
for $\lambda \in \R$. Since further $f_{Y_aY_C}(\lambda) f_{Y_CY_C}(\lambda)^{-1} f_{Y_CY_a}(\lambda) \geq 0$ and the integral is monotonous, we obtain the integrability
\begin{align*}
& \uint \left\vert f_{Y_a Y_C}(\lambda) f_{Y_CY_C}(\lambda)^{-1}  f_{Y_CY_C}(\lambda) \overline{f_{Y_a Y_C}(\lambda) f_{Y_CY_C}(\lambda)^{-1} }^\top \right\vert d\lambda \nonumber \\
&\quad = \uint f_{Y_aY_C}(\lambda) f_{Y_CY_C}(\lambda)^{-1} f_{Y_CY_a}(\lambda) d\lambda  \\
&\quad \leq \uint f_{Y_aY_a}(\lambda) d\lambda <\infty,
\end{align*}
where the finiteness follows from \Cref{properties of spectral density}(a). In summary, $\widehat{Y}_{a\vert C}(t) \in \mathcal{L}_{Y_C}$ for $t\in \R$.\\
\textsl{Step 2: \:} Due to \cite{RO67}, I, (7.2), any element $Y^C \in \mathcal{L}_{Y_C}$ has the spectral representation
\begin{align*}
Y^C = \uint \varphi(\lambda) \Phi_C(d\lambda) \quad \text{$\mathbb{P}$-a.s.,}
\end{align*}
 where $\varphi\in L^2\left( f_{Y_C Y_C} \right)$. Now, writing $Y_a(t)$ in its spectral representation \eqref{spectral representation of stationary process}, it holds that
 \begin{align*}
     &\BE \bigg[ \Big( Y_a(t) - \widehat{Y}_{a\vert C}(t) \Big) \overline{Y^C} \bigg] \\
     &\quad = \BE \left[ \left( \uint e^{i\lambda t} \Phi_a(d\lambda) - \uint e^{i \lambda t} f_{Y_aY_C}(\lambda) f_{Y_CY_C}(\lambda)^{-1} \Phi_C(d \lambda) \right) \overline{\uint \varphi(\lambda) \Phi_C(d\lambda)}
     \right] \\
     &\quad= \uint e^{i\lambda t}  f_{Y_a Y_C}(\lambda) \overline{\varphi(\lambda)}^\top d\lambda
     - \uint e^{i \lambda t} f_{Y_a Y_C}(\lambda) f_{Y_CY_C}(\lambda)^{-1}
     f_{Y_CY_C}(\lambda) \overline{\varphi(\lambda)}^\top d\lambda
=0.
 \end{align*}
Thus, $Y_a(t)- \widehat{Y}_{a\vert C}(t) \in \mathcal{L}_{Y_C}^\perp$ for $t\in \R$.\\
\textsl{Step 3: \:} Since \mbox{$\mathcal{L}_{Y_C} = \mathcal{L}_{Y_C}^*$} (\Cref{Gleichheit linearer Raeume}) the optimisation problem \eqref{optimisation problem} is equivalent to
\begin{align*}
    \underset{Y^C\in \mathcal{L}_{Y_C}}{\text{min}} \BE \left[ \left\vert Y_a(t) - Y^C \right\vert^2 \right].
\end{align*}
From the minimality property of the orthogonal projection, we obtain \mbox{$Y^C=P_{\mathcal{L}_C} Y_a(t)$} is the optimal solution to this optimisation problem. Due to Step 1 and Step 2 the function $\varphi_{a \vert C}(\lambda)=f_{Y_a Y_C}(\lambda) f_{Y_CY_C}(\lambda)^{-1}$, $\lambda \in \R$, is then the optimal function in \eqref{optimisation problem}.
\end{proof}

\begin{proof}[Proof of \Cref{remainder satisfies Assumption 1}]
We can write
\begin{align*}
\ovA(t)
&= \uint e^{i \lambda t} \Phi_A(d\lambda) - \uint e^{i \lambda t} f_{Y_AY_C}(\lambda) f_{Y_CY_C}(\lambda)^{-1} \Phi_C(d\lambda) \\
&= \uint e^{i \lambda t} \left( E_A^\top - f_{Y_AY_C}(\lambda) f_{Y_CY_C}(\lambda)^{-1} E_C^\top \right) \Phi_V(d\lambda),
\end{align*}
where $E_A \in \R^{k \times |A|}$ (and analogously $E_C \in \R^{k \times |C|}$) is the matrix defined by itsentries
    \begin{align*}
    [E_A]_{ij} =
    \begin{cases}
    1, & i=j, \: i,j \in A,\\
    0, & \text{else}.
    \end{cases}
    \end{align*}
Therefore, the noise process $(\ovA(t))_{t\in \R}$ is a linear transformation of the wide-sense stationary process $\CY_V$ with spectral characteristic $E_A^\top - f_{Y_AY_C}(\lambda) f_{Y_CY_C}(\lambda)^{-1} E_C^\top$, $\lambda\in \R$. Thus $\CY_V$ is also wide-sense stationary \cite[I, (8.2)]{RO67}. Furthermore, the linear transformation has a spectral density function, which is given by \cite[I, (8.13)]{RO67}
\begin{align*}
\fAA(\lambda)
& = \left(E_A^\top - f_{Y_AY_C}(\lambda) f_{Y_CY_C}(\lambda)^{-1} E_C^\top \right) f_{Y_V Y_V}(\lambda) \overline{\left(E_A^\top - f_{Y_AY_C}(\lambda) f_{Y_CY_C}(\lambda)^{-1} E_C^\top \right)}^\top \\
& = f_{Y_A Y_A}(\lambda)
- f_{Y_A Y_C}(\lambda) f_{Y_C Y_C}(\lambda)^{-1} f_{Y_C Y_A}(\lambda).
\end{align*}
Then \Cref{properties of spectral density}(b) yields
\begin{align*}
\cAA(t) = \uint e^{i \lambda t} \left(f_{Y_A Y_A}(\lambda) - f_{Y_A Y_C}(\lambda) f_{Y_CY_C}(\lambda)^{-1}  f_{Y_C Y_A}(\lambda)  \right) d\lambda
\end{align*}
for $t\in \R$. In particular, the spectral density function of $(\ovAB(t))_{t\in \R}$ is given by
\begin{align*}
 \fABAB(\lambda)
 = f_{Y_{A \cup B}Y_{A \cup B}}(\lambda) - f_{Y_{A \cup B}Y_{C}}(\lambda) f_{Y_{C}Y_{C}}(\lambda)^{-1} f_{Y_{C}Y_{A \cup B}} (\lambda).
\end{align*}
Thus, the cross-spectral density function is, for almost all $\lambda \in \R$,
\begin{align*}
\fAB(\lambda)
& = E_A^\top \left( f_{Y_{A \cup B}Y_{A \cup B}}(\lambda) - f_{Y_{A \cup B}Y_{C}}(\lambda) f_{Y_{C}Y_{C}}(\lambda)^{-1} f_{Y_{C}Y_{A \cup B}} (\lambda) \right) E_B \\
& = f_{Y_A Y_B}(\lambda)
- f_{Y_A Y_C}(\lambda) f_{Y_C Y_C}(\lambda)^{-1} f_{Y_C Y_B}(\lambda),
\end{align*}
and, for all $t\in \R$, it holds that
\begin{equation*}
    \cAB(t) = \uint e^{i \lambda t} \left(f_{Y_A Y_B}(\lambda) - f_{Y_A Y_C}(\lambda) f_{Y_CY_C}(\lambda)^{-1}  f_{Y_C Y_B}(\lambda)  \right) d\lambda. \qedhere
\end{equation*}
\end{proof}

\begin{proof}[Proof of \Cref{characterisation with spectral density function}]
Suppose that $\CY_A \inde  \CY_B \: \vert \:  \CY_{C}$. By definition of this relation we obtain the first characterisation $\cAB(t) =0_{|A| \times |B|}$ for $t \in \R$. For the second characterisation suppose that \mbox{$\cAB(t)=0_{|A| \times |B|}$} for $t \in \R$. Then the Fourier inversion formula \cite[Proposition 2.2.37]{PI09} yields \mbox{$\fAB(\lambda) = 0_{|A| \times |B|}$} for almost all $\lambda \in \R$. If $\fAB(\lambda) = 0_{|A| \times |B|}$ for almost all $\lambda \in \R$, then \Cref{properties of spectral density} gives $ \cAB(t) = 0_{|A| \times |B|}$ for $t \in \R$.
For the third characterisation, suppose that $\fAB(\lambda) = 0_{|A| \times |B|}$ for almost all $\lambda \in \R$. Then $\RAB(\lambda) = 0_{|A| \times |B|}$ holds by \Cref{partial correlation matrix}. If we additionally assume that $\fAA(\lambda)>0$ and $\fBB(\lambda)>0$, then \Cref{partial correlation matrix} provides the second direction.
\end{proof}

\begin{proof}[Proof of \Cref{graphoid}]\mbox{}\\
\textsl{(P4)}\: The relations $\CY_A \inde \CY_B \: \vert \: \CY_D$, $\CY_A \inde \CY_C \: \vert \: (\CY_B,\CY_D)$, and \Cref{characterisation with inverse}
result in
\begin{align}\label{eqn 3.1}
\left[g_{Y_{A \cup B \cup D}Y_{A \cup B \cup D}}(\lambda)\right]_{A B}=0_{|A| \times |B|} \quad \text{and} \quad
\left[g_{Y_{A \cup B \cup C \cup D}Y_{A \cup B \cup C \cup D}}(\lambda)\right]_{A C}=0_{|A| \times |C|}
\end{align}
for almost all $\lambda \in \R$. Along with \Cref{deleting confounder}, we obtain
\begin{align}\label{eqn 3.2}
0_{|A| \times |B|}
=&\: \left[g_{Y_{A \cup B \cup D}Y_{A \cup B \cup D}}(\lambda)\right]_{A B} \nonumber \\
=&\:\left[g_{Y_{A \cup B \cup C \cup D}Y_{A \cup B \cup C \cup D}}(\lambda)\right]_{A B} -
 \left[g_{Y_{A \cup B \cup C \cup D}Y_{A \cup B \cup C \cup D}}(\lambda)\right]_{A C} \nonumber \\
 &\:\left[\left(g_{Y_{A \cup B \cup C \cup D}Y_{A \cup B \cup C \cup D}}(\lambda)\right]_{C C}\right)^{-1}
 \left[g_{Y_{A \cup B \cup C \cup D}Y_{A \cup B \cup C \cup D}}(\lambda)\right]_{C B} \nonumber \\
=&\:\left[g_{Y_{A \cup B \cup C \cup D}Y_{A \cup B \cup C \cup D}}(\lambda)\right]_{A B}
\end{align}
for almost all $\lambda \in \R$. In summary, equations \eqref{eqn 3.1} and \eqref{eqn 3.2} give
\begin{align*}
\left[g_{Y_{A \cup B \cup C \cup D}Y_{A \cup B \cup C \cup D}}(\lambda)\right]_{A (B \cup C)}
=0_{|A| \times (|B| + |C|)}
\end{align*}
for almost all $\lambda \in \R$. \Cref{characterisation with inverse} implies $\CY_A \inde (\CY_B,\CY_C) \: \vert \: \CY_D$.
\end{proof}

\subsection{Proofs of Section~\ref{sec: Partial correlation graphs}}\label{sec: Proofs of sec: Partial correlation graphs}

\begin{proof}[Proof of \Cref{Zusammenhang zu causality graph}]
 Theorem 5.15 in \cite{VF23pre} provides that $\{a\} \msep \{b\} \: \vert\: V\setminus\{a,b\} \:\:\: [G_{GC}]$ implies $ \mathcal{L}_{Y_a} \perp \mathcal{L}_{Y_b} \: \vert \: \mathcal{L}_{Y_{V\setminus\{a,b\}}}$. This conditional orthogonality relation immediately implies $L_{Y_a}(t) \perp L_{Y_b}(t) \: \vert \: \mathcal{L}_{Y_{V\setminus\{a,b\}}}$ for all $t\in \R$ by subset arguments, which in turn yields $a\edge b \notin E_{PC}$ due to \Cref{Bedingte Orthogonalitaet}.
\end{proof}

\begin{proof}[Proof of \Cref{Zusammenhang 2 zu causality graph}]
By definition and due to \cite{EI11}, Theorem 3.1 and Lemma 3.2, we obtain that
\begin{align*}
a\edge b \notin E^a_{GC}
     \quad &\Leftrightarrow \quad
    \text{$a$ and $b$ are not collider connected in $G_{GC}$}, \\
    \quad &\Leftrightarrow \quad
    \dis\left( a \cup \ch(a)\right)\cap \dis\left( b \cup \ch(b) \right) \: \text{in $G_{GC}$},\\
    \quad &\Leftrightarrow \quad 
    \{a\} \msep \{b\} \: \vert\: V\setminus\{a,b\} \:\:\: [G_{GC}], \\
    \quad &\Leftrightarrow \quad 
    \{a\} \bowtie \{b\} \: \vert\: V\setminus\{a,b\} \:\:\: [G_{GC}^a].
\end{align*}
These statements are then of course all sufficient for $a\edge b \notin E_{PC}$ due to the previous \Cref{Zusammenhang zu causality graph}, and $E_{PC} \subseteq E^{a}_{GC}$ is valid.
\end{proof}

\subsection{Proofs of Section~\ref{sec: Partial correlation graphs for MCAR processes}}\label{sec: Proofs of sec: Partial correlation graphs for MCAR processes}

\begin{proof}[Proof of \Cref{Characterisierung PC fuer MCAR}]
First of all, we insert the AR polynomial $P$ in $g_{Y_VY_V}(\lambda)$ to get
\begin{align*}
g_{Y_VY_V}(\lambda)
=& \: 2\pi \left(\sum_{m=0}^{p} A_{p-m}^\top (-i\lambda)^m \right) \BS_L^{-1} \left(\sum_{\ell=0}^p A_{p-\ell} (i\lambda)^\ell \right) \\
=& \: 2\pi \sum_{n=0}^{2p} \sum_{m=0\vee n-p}^{n\wedge p} (-1)^m A^\top _{p-m}\BS_L^{-1}A_{p-n+m} (i\lambda)^n.
\end{align*}
In the last step, we arrange the addends according to the degree of $\lambda$ and substitute $n=\ell+m$, where $n=0,\ldots,2p$. Since $0 \leq \ell = n-m \leq p$ and $0 \leq m \leq p$, we obtain the boundary $0\vee n-p \leq m \leq n\wedge p$. Since the components of $g_{Y_VY_V}$ are polynomials, the components are zero functions if and only if the corresponding coefficients are zero. Then, by \Cref{Prop: Definition partial correlation graph for MCAR}, we obtain that
\begin{align}\label{Charakterisierung allgemein}
a\edge b \notin E_{PC}
\:\: \Leftrightarrow \:\:
\left[ \sum_{m=0\vee n-p}^{n\wedge p} (-1)^m A^\top _{p-m}\BS_L^{-1}A_{p-n+m}
\right]_{ab}=0 \:\: \text{for $n=0,\ldots,2p$.}
\end{align}
\textsl{(i)\:} Assume that $\BS_L = \sigma^2 I_k$. Then $\BS_L^{-1} = 1/\sigma^2 I_k$ holds and since $\sigma^2>0$, relation \eqref{Charakterisierung allgemein} is equivalent to \Cref{Characterisierung PC fuer MCAR}(i).\\
\textsl{(ii)\:} Assume that $A_m$, $m=1,\ldots,p$, are diagonal matrices. Then the AR polynomial $P$ is a diagonal matrix polynomial and $a\edge b \notin E_{PC}$ is equivalent to
\begin{align*}
0 = \left[P(-i\lambda) \BS_L^{-1} P(i\lambda)\right]_{ab}
&= \left[P(-i\lambda)\right]_{aa} \left[\BS_L^{-1}\right] _{ab} \left[P(i\lambda)\right]_{bb}
\end{align*}
for all $\lambda \in \R$. Due to the causality assumption $\sigma(\BA) \subseteq (-\infty,0) +i\R$ and the structure of $\BA$, the diagonal matrix $A_p$ is not singular and in particular the diagonal elements of $A_p$ are not zero. Thus the diagonal elements of $P(i\lambda)$ are never zero and $a\edge b \notin E_{PC}$ is equivalent to $[\BS_L^{-1}] _{ab}=0$.
\end{proof}

\begin{proof}[Proof of \Cref{Spezialfall Beziehung}]
The assumptions that $\Sigma_L=I_k$, $\BA$ is symmetric, and \eqref{eq5.3} imply
\begin{gather*}
\left[ \BA\right]_{ba} - \left[\BA \right]_{ab}
=0, \quad
\left[ \BA^\top \BA \right]_{ab}
= \left[ \BA \BA^\top \right]_{ab}
=0.
\end{gather*}
Thus, \eqref{OU Kriterium 1b} yields $a \edge b \notin E_{PC}$.
\end{proof}

\begin{proof}[Proof of \Cref{Separationsvergleich local}]
The equivalences
\begin{align*}
 a\edge b \notin E^{0,a}_{GC}
     \quad &\Leftrightarrow \quad
    \text{$a$ and $b$ are not collider connected in $G_{GC}^0$},\\
    \quad &\Leftrightarrow \quad
    \dis\left( a \cup \ch(a)\right)\cap \dis\left( b \cup \ch(b) \right) \: \text{in $G_{GC}^0$},\\
    \quad &\Leftrightarrow \quad
    \{a\} \bowtie \{b\} \: \vert\: V\setminus\{a,b\} \:\:\: [G_{GC}^{0,a}], \\
    \quad &\Leftrightarrow \quad
    \{a\} \msep \{b\} \: \vert\: V\setminus\{a,b\} \:\:\: [G_{GC}^{0}],
\end{align*}
were already established in \Cref{Zusammenhang 2 zu causality graph}, regardless of the specific definition of the graphical model. Thus, only $a\edge b \notin E_{PC}$ needs to be proved, which we do by contradiction. Suppose that $a\edge b \in E_{PC}$. Then there exists a $\lambda\in \R$, such that
\begin{equation*}
0 \neq \left[P(-i\lambda)^\top \Sigma_L^{-1} P(i\lambda) \right]_{ab}
  = \sum_{c\in V} \sum_{d\in V} \left[P(-i\lambda)\right]_{ca} \left[ \Sigma_L^{-1} \right]_{cd} \left[P(i\lambda)\right]_{db}.
\end{equation*}
Consequently, there exist vertices $c,d\in V$, such that
\begin{align*}
\left[P(-i\lambda)\right]_{ca}\neq 0, \quad
\left[ \Sigma_L^{-1} \right]_{cd} \neq 0, \quad
\left[P(i\lambda)\right]_{db} \neq 0.
\end{align*}
This means that there are directed edges $a\rarrow c$ and $b\rarrow d$ in $G_{GC}^0$. 

{If $c=b$ (or $d=a$), then the edge $a\rarrow b$ ($b\rarrow a$) is trivially a collider path which is contradiction and hence, $a\edge b \notin E_{PC}$ holds. 

Thus, in the following we assume that $c\not=b$ and $d\not=a$.
Then} there exists a path $\pi$ between $c$ and $d$ of only undirected edges in the local causality graph $G_{GC}^0$ \citep[p.~341]{EI07}. Indeed, for an MCAR(1) process with $\BA=-I_k$ that is driven by the same Lévy process $L$, we have $ a \edge b \notin \widetilde{E}_{PC}$ if and only if $[\Sigma_L^{-1}]_{ab}=0$ and $a \inst b \notin  \widetilde{E}_{GC}$ if and only if $[\Sigma_L]_{ab}=0$. Additionally, there are no directed edges in the causality graph $\widetilde{G}_{GC}$. Then a consequence of \Cref{Zusammenhang 2 zu causality graph} is that $a\edge b \in \widetilde{E}_{PC}$ ($[\Sigma_L^{-1}]_{cd}\neq 0$) implies $a\edge b \in \widetilde{E}_{GC}^a$ and $\dis(a)\cap \dis(b) \neq \emptyset$ in $\widetilde{E}_{GC}$.
Thus, there exists a path $\pi$ of only undirected edges between $c$ and $d$ in the causality graph $\widetilde{G}_{GC}$, i.e., for some $c=\alpha_1,\ldots,\alpha_l=d\in V$ we have $[\Sigma]_{\alpha_i\alpha_{i+1}}\not=0$ for $i=1,\ldots,l-1$ 
and hence, in the local causality graph ${G}_{GC}^0$.

We complete $\pi$ with the directed edges to get a path $\widetilde{\pi}$ between $a$ and $b$ on which every intermediate vertex is a collider. This is a contradiction of the premise and the statement $a\edge b \notin E_{PC}$ holds.
\end{proof}

\end{document}